\documentclass[12pt,a4paper,leqno]{article}

\usepackage{amssymb}
\usepackage{german}

\def\mwmatrix#1{\null\,\vcenter{\openup1\jot \mathsurround00pt\ialign{&\strut\hfil$\displaystyle{{}##{}}$\hfil\crcr#1\crcr}}\,}

\newtheorem{Proposition}{Proposition}[section]
\newtheorem{Definition}[Proposition]{Definition}
\newtheorem{Satz}[Proposition]{Satz}
\newtheorem{Lemma}[Proposition]{Lemma}
\newtheorem{Korollar}[Proposition]{Korollar}
\newtheorem{Bemerkung}[Proposition]{Bemerkung}
\newtheorem{Notation}[Proposition]{Notation}
\newtheorem{Beispiel}[Proposition]{Beispiel}

\newcommand {\F}{\mbox{${\cal F}$}}

\newcommand {\Abb}[5]{\begin{eqnarray*} #1  #2  & \longrightarrow & #3 \\ 
              #4 & \longmapsto & #5 \end{eqnarray*} }
\newcommand {\beweis}{\noindent {\bf Beweis: }}

\newcommand {\aequi}{\mbox{$\Leftrightarrow$}}
\newcommand {\Fl}{\mbox{${\cal F} \!\!\!\!\!\!\!\; {\cal J}$}}
\newcommand {\su}{\mbox{$ \subset $}}
\newcommand {\nsu}{\mbox{$ \not\subset $}}
\newcommand{\B}{{\cal B}}
\newcommand{\usupp}[1]{\underline{\mbox{supp}}(#1)}
\newcommand{\supp}[1]{\mbox{supp}(#1)}
\newcommand{\subn}{\mathop{\subset}\limits_{\neq}}
\newcommand {\exseq}[3]{\mbox{$0 \longrightarrow #1 \longrightarrow #2
    \longrightarrow #3 \longrightarrow 0 $ } }
\newcommand {\maxwedge}[1]{\mbox{$\bigwedge^{\max} #1 $}}
\newcommand {\slo}{\mbox{slope}}
\newcommand {\qed}{\hfill $\square$}

\newenvironment{rnumer}{\begin{list}{\roman{enumi})}{\usecounter{enumi}\setlength{\leftmargin}{0.5cm}}}{\end{list}}

\begin{document}

\title{Kohomologie von Periodenbereichen "uber endlichen K"orpern}
\author{Sascha Orlik}
\date{}

\maketitle

\begin{abstract} Periodenbereiche sind gewisse offene Unterr"aume von
  verallgemeinerten Flaggenvariet"aten, welche durch
  Semistabilit"atsbedingungen beschrieben werden. In dem Fall eines
  endlichen Grundk"orpers bilden diese eine Zariski-offene Untervariet"at,
  im Fall eines lokalen K"orpers einen zul"assigen offenen Unterraum im
  Sinne der rigiden algebraischen Geometrie. In dieser Arbeit berechnen wir
  f"ur den  Fall eines endlichen Grundk"orpers die $\ell$-adische
  Kohomologie  mit kompaktem Tr"ager dieser  Periodenbereiche .
\end{abstract}


\centerline{\large \bf Einleitung}

\bigskip
Der Begriff des Periodenbereiches wurde von Griffiths \cite{G} eingef"uhrt.
Diese Objekte stellen gewisse offene Untermengen von
verallgemeinerten Flaggenvariet"aten "uber $\mathbb C$  dar,
welche polarisierte $\mathbb R$-Hodgestruk\-turen eines gegebenen Typs
parametrisieren. Hierdurch motiviert hat Rapoport in \cite{R1}
Periodenbereiche "uber $p$-adischen beziehungsweise endlichen K"orpern
definiert. Wir wollen die letztere Variante kurz erl"autern. Sei dazu $k$ ein endlicher
K"orper. Ein filtrierter $k$-Vektorraum ist ein Paar $(V,\F^\bullet)$ bestehend aus
einem endlich-dimensionalen $k$-Vektorraum $V$ und einer $\mathbb R$-Filtration
$\F^\bullet$ auf $V,$ welche "uber einer endlichen K"orpererweiterung von $k$ definiert
ist. Analog zu dem Fall von Vektorb"undeln auf einer Kurve existiert der Begriff des Anstieges. F"ur einen filtrierten Vektorraum
$V=(V,\F^\bullet)\neq 0$ ist dieser definiert durch
$$\slo(V)=\slo_{\F^\bullet}(V)=\frac{\sum_x x\dim gr_{\F^\bullet}^x(V)}{\dim V}.$$ 
Entsprechend hei"st eine Filtration semistabil, wenn f"ur jeden nicht trivialen $k$-rationalen
Unterraum $U$ die Ungleichung  $\slo(U) \leq \slo(V)$ erf"ullt ist. Dabei ist
der Unterraum $U$ mit der induzierten Filtration versehen. Man betrachtet nun
numerierte Filtrationen eines gegebenen Typs auf einem $k$-Vektorraum $V$,
d.h. man fixiert eine Funktion $g:\mathbb R \rightarrow \mathbb Z_{\geq 0}$
mit endlichem Tr"ager und l"a"st nur solche Filtrationen $\F^\bullet$ auf $V$ zu,
f"ur die $\dim gr^x_{\F^\bullet}(V)=g(x)\; \forall x\in\mathbb R .$  
Diese lassen sich durch eine "uber $k$ definierte
Flaggenvariet"at $\Fl_g$ parametrisieren. Der Unterraum der semistabilen
Filtrationen hei"st dann der zu $g$ geh"orige Periodenbereich und wird
mit $\Fl_g^{ss}$ notiert \cite{R1}. Bei diesem Raum handelt es sich um
einen  offenen Teilraum von  $\Fl_g.$  Somit ist er mit der Struktur
einer "uber $k$ definierten Variet"at versehen. In der Situation eines
$p$-adischen Grundk"orpers $k$ betrachtet man anstatt filtrierter
Vektor"aume filtrierte Isokristalle "uber $k$, f"ur die analog der Begiff
des Anstieges existiert. Hier entsprechen die
semistabilen Objekte vom Anstieg 0  genau den schwach zul"assigen Isokristallen
(vgl. \cite{RZ}) aus der Fontainschen Theorie (vgl. \cite{Fo}). Die Theorie der
Periodenbereiche  l"a"st  sich sowohl f"ur endliche als auch f"ur $p$-adische
Grundk"orper auf den Fall von beliebigen reduktiven Gruppen verallgemeinern.

Das prominenteste Beispiel eines Periodenbereiches ist der sogenannte
Drinfeldsche Halbraum $\Omega^d$ der Dimension $d-1.$ Dieser ist im Falle eines endlichen
Grundk"orpers $k$ das Komplement aller $k$-rationalen Hyperebenen im
projektiven Raum $\mathbb P^{d-1}_k,$ d.h. es ist 
$$ \Omega^d=\mathbb P_k^{d-1}\setminus \bigcup_{H \subset k^{d} \atop  \dim H=d-1}\mathbb P(H).$$
Dieser Raum ist im Gegensatz
zu den meisten Periodenbereichen sogar affin. In der $p$-adischen Situation hat
man bei diesem Beispiel lediglich den Grundk"orper $k$ durch den entsprechenden
lokalen K"orper zu ersetzen.
Eine ausf"uhrliche Beschreibung dieser $p$-adischen Variante des
Drinfeldschen Halbraumes wird in der Arbeit \cite{SS} diskutiert. 
Es sei schlie"slich erw"ahnt, da"s es wohl keinen allgemeinen Zusammenhang der
Periodenbereiche zu den Deligne-Lusztig-Variet"aten \cite{DL} gibt, obgleich
man den Drinfeldschen Halbraum als eine solche beschreiben kann
(vgl. \cite{DL} 2.2).

Ein nat"urliches Anliegen ist es, die Kohomologie der Periodenbereiche f"ur
beliebige reduktive Gruppen $G$ zu kennen. Dabei verstehen wir unter der
Kohomologie die $\ell$-adische Kohomologie mit kompaktem Tr"ager. Die
Periodenbereiche sind im allgemeinen Fall erst "uber einer endlichen
K"orpererweiterung $E$ des Grundk"orpers definiert. 
Die Kohomologiegruppen sind dann auf
nat"urliche Weise mit Operationen sowohl von der Gruppe $G(k)$, f"ur $k$
endlich, (bzw. $J(k)$, im Fall eines lokalen K"opers $k$)  als auch der
Galoisgruppe $Gal(\overline{E}/E)$  
versehen. Dabei bezeichnet  $J$ eine gewisse innere Form der Gruppe $G.$
In der $p$-adischen Variante wurde die
$\ell$-adische Kohomologie des Drinfeldschen Halbraumes in  \cite{SS}
berechnet.  Es gilt 
$$H^{\ast}_c(\Omega^d,\mathbb Q_\ell)=\bigoplus_{i=0}^{d-1}
v^G_{P_{(d-i,1,\ldots,1)}}(i+1-d)[-2(d-1)+i],$$
wobei $P_{(d-i,1,\ldots,1)}$ die der Partition $(d-i,1,\ldots,1)$
entsprechende standard-parabolische Untergruppe der $GL_d$ ist und
 $v^G_P,$ f"ur $P\subset G$ parabolisch,  die
Darstellung $Ind^{G(k)}_{P(k)}\mathbb Q_{\ell}/\sum_{P \subn
  Q}Ind^{G(k)}_{Q(k)}\mathbb Q_{\ell}$ bezeichnet (vgl. \cite{SS}). 
Die Notation $[-n], n\in \mathbb N,$ bedeutet hierbei, da"s
der voranstehende Modul in den Grad $n$ des Kohomologieringes geshiftet
wird. Unter der Bezeichnung $(m), m\in \mathbb Z,$ verstehen wir den
$m$-fachen Tate-Twist. Die Formel ist ebenso im Fall eines endlichen
Grundk"orpers g"ultig. 

F"ur beliebige Periodenbereiche haben Kottwitz und
Rapoport Formeln f"ur die Euler-Poincar\'e-Charakteristik bzgl. der
Grothendieck-Gruppe von $G(k)\times Gal(\overline{E}/E)$- (bzw. $J(k)\times Gal(\overline{E}/E)$-) Darstellungen 
hergeleitet (vgl. \cite{R1}-\cite{R3}). Das Hauptresultat dieser Arbeit verallgemeinert die obige Formel im Fall
$G=GL(V)$ f"ur endliche Grundk"orper $k$. Sie best"atigt eine Vermutung von Kottwitz und Rapoport. Der
Beweis ist allerdings v"ollig anders als der von Rapoport und Kottwitz f"ur
die Euler-Poincar\'e-Formel. 

Wir wollen das Ergebnis nun erl"autern. In unserem  Fall
ist $E=k.$ Bezeichne mit $W$ die Weylgruppe von $G$ bez"uglich des
Diagonaltorus $T.$ Es sei $\mu=\mu_g=(x_1^{g(x_1)},\ldots,x_r^{g(x_r)}) \in
X_\ast(T)_{\mathbb   R}$ der zu $g$ assoziierte reelle Cocharakter, wobei
supp($g$)=$\{x_1> \ldots >x_r\}$ den Tr"ager von $g$ darstellt. Ferner sei
$W_\mu$ der Stabilisator von $\mu$ in $W$ und $W^\mu$ die Menge
der Kostant-Repr"asentanten von $W/W_\mu.$ Sei $\Delta$ die zur Borelgruppe
$B$ der oberen Dreiecksmatrizen geh"orige Wurzelbasis. Seien
$\omega_\alpha, \alpha \in \Delta,$ die assoziierten Fundamentalgewichte. Ordnet man jedem $w \in
W^\mu$ die Teilmenge $\Delta_w=\{\alpha \in \Delta; \langle w\mu,
\omega_\alpha\rangle\;>0\}$ von $\Delta$ zu, so lautet unsere Formel f"ur die
Kohomologie:

\bigskip
\noindent {\bf Satz:} $\;\;\;\;\;\;\;H^\ast_c(\Fl_g^{ss}, \mathbb
Q_\ell)=\bigoplus_{w \in W^\mu} v^G_{P_w}(-l(w))[-2l(w)-\#\Delta_w].$

\noindent Hierbei notieren wir mit  $P_w$ diejenige
standard-parabolische Untergruppe, deren unipotentes Radikal von $\Delta_w$ erzeugt wird.  F"ur ein $w\in W$ mit
  $\Delta_w=\Delta$ erhalten wir somit die Steinberg-Darstellung $v^G_B.$
  Anders als die Steinberg-Darstellung sind die $v^G_{P_w}$ jedoch im
  allgemeinen nicht mehr irreduzibel. Es ist zu vermuten, da"s die obige Formel auch f"ur beliebige reduktive Gruppen
und auch f"ur die $p$-adische Situation gilt. Wie man der Formel entnimmt, ist die
  Kohomologie der Periodenbereiche im allgemeinen nicht mehr rein. 
Jedes $w \in W^\mu$ induziert einen Beitrag zur Kohomologie,
welcher sich im Grad $2l(w)+\# \Delta_w$ des Kohomologieringes
befindet. Dabei verhalten sich diese beiden Summanden bez"uglich der
Bruhatordnung auf $W$ gegenl"aufig, d.h. aus $w' \leq w$ folgt $l(w') \leq l(w),$ aber
$\Delta_w \subset \Delta_{w'}.$ Deswegen kann es, anders als
bei der Kohomologie von Flaggenvariet"aten, durchaus passieren, da"s f"ur
$w' \leq w$ der induzierte Grad von $w'$ gr"o"ser als derjenige
von $w$ ist. 
Trotzdem l"a"st sich, wie Rapoport bemerkt hat, aus der Formel ein Verschwindungssatz ableiten.

{\noindent \bf Korollar:} {\it  F"ur den zu $g$ assoziierten Periodenbereich $\Fl^{ss}_g$ ist
$$ H^i_c(\Fl_g^{ss},\mathbb Q_\ell)=0 ,\;\;\; 0 \leq i \leq d-2$$
und}
$$ H^{d-1}_c(\Fl_g^{ss},\mathbb Q_\ell)= v^G_B .$$

Wir wollen  nun kurz die Beweisstrategie des obigen Resultats 
schildern. Eine ausf"uhrlichere Beschreibung findet man im Paragraphen 4. Es ist etwas "uberraschend, da"s die Berechnung der Kohomologie
mit der Methode dieser Arbeit gegl"uckt ist. Wir schreiben zun"achst
die Menge der instabilen Punkte in der Flaggenvarit"at als endliche
Vereinigung von abgeschlossenen Untervariet"aten $Y_U,$ die durch die Menge
der $k$-rationalen Unterr"aume $U \subset V$ parametrisiert werden. Dabei besteht $Y_U$
aus denjenigen Filtrationen, f"ur die der Unterraum $U$ die
Semistabilit"atsungleichung verletzt. Der n"achste nat"urliche Schritt
w"are die Beschreibung s"amtlicher Durchschnitte obiger
Untervariet"aten, um anschlie"send die Kohomologie mittels der
verallgemeinerten Mayer-Vietoris Sequenz zu berechnen. Jedoch ist dies f"ur
allgemeines $d$ ein sehr schwieriges Problem. Es ist noch nicht einmal bekannt,
f"ur welche $U_1,\ldots, U_r$ der Durchschnitt $Y_{U_1}\cap \ldots \cap
Y_{U_r}$ nicht leer ist. Jedoch lassen
sich gewisse Durchschnitte beschreiben, mit deren Hilfe  man
"uberraschenderweise die Kohomolgie von $Y$ und damit die des offenen
Komplementes  berechnen kann. Der Punkt ist, da"s f"ur ein fixiertes
Element $x$  aus der Flaggenvariet"at die Menge der $U,$ f"ur die $x$ in
$Y_U$ enthalten ist, einen kontrahierbaren Unterkomplex des Tits-Komplexes
zu $GL(V)(k)$ induziert. 
Dieses Ph"anomen entspricht dem Ergebnis von Mumford "uber die
Konvexit"at der Menge von Einparameteruntergruppen, welche die
Semistabilit"at verletzen (vgl. \cite{M} 2.2).

Wir kommen nun zur Inhaltsangabe dieser Arbeit.
Der Paragraph 1 behandelt die Theorie der numerierten Filtrationen. 
Wir beschr"anken uns hierbei auf die im weiteren wichtigen Definitionen und
Eigenschaften. Im Paragraphen 2 f"uhren wir Variet"aten $\Fl_g(\B')$ ein,
deren Kohomologie wir berechnen wollen. Dabei ist $\B'$ eine Teilmenge
einer Menge $\B$ von sogenannten Unterfunktionen von $g$. Das
abgeschlossene Komplement $Y(\B')$ von $\Fl_g(\B')$ ist die Vereinigung von
konstruierbaren Untermengen $Y(h), \; h\in \B'.$ Auf der Menge
$\B$ definieren wir eine Ordnungsrelation, mit der wir den
Abschlu"s der $Y(h)$ in $\Fl_g$ beschreiben. 
Anschlie"send zeigen wir, da"s die Periodenbereiche  Spezialf"alle
der Variet"aten $\Fl_g(\B')$ sind. 
Im Paragraphen 3 zeigen wir, da"s die Ordnung auf $\B$
zu der Bruhatordnung der Weylgruppe von $GL(V)$ korres\-pondiert. 
Schlie"slich beweisen wir die obige
Abschlu"seigenschaft der Mengen $Y(h)$, indem wir uns auf den Fall von
Schubertzellen zur"uckziehen. Die Hauptresultate der Arbeit werden im 
Paragraph 4
formuliert. Daraus leiten wir den bereits erw"ahnten Verschwindungssatz ab. 
Anschlie"send diskutieren wir einige Eigenschaften der Kohomologie von
Periodenbereichen, wie zum Beispiel die Reinheit.
Die Konstruktion des fundamentalen Komplexes von \'etalen Garben auf
$Y(\B')$ findet man im Paragraphen 5. Die Azyklizit"at dieses Komplexes
beweisen wir im Paragraphen 6. Im letzten Paragraphen werten wir schlie"slich die
Spektralsequenz aus, welche durch den azyklischen Komplex induziert
wird. Diese Auswertung liefert die gesuchte Kohomologie der Variet"aten
$Y(\B').$ 

An dieser Stelle m"ochte ich erw"ahnen, da"s es sich bei dieser Arbeit um eine
Dissertation der Mathematisch-Naturwissenschaftlichen Fakult"at der
Universit"at zu K"oln handelt. Ich m"ochte  mich recht herzlich bei
Prof. Dr. M. Rapoport f"ur die Vergabe und intensive Betreuung dieser Arbeit bedanken. Au"serdem m"ochte ich mich bei T. Fimmel, T. Wedhorn und
O. B"ultel und U. G"ortz bedanken, mit denen ich rege Diskussionen f"uhrte
und die mir sehr wertvolle Tips gaben.

\section{Numerierte Filtrationen}
Wir erinnern kurz an einige Begriffe aus der Theorie der numerierten
Filtrationen. Als Referenz seien dabei die Arbeiten  \cite{R1}, \cite{R2},
\cite{FW} und \cite{O}  erw"ahnt. 

In diesem Paragraphen bezeichnen wir mit $k$ einen beliebigen
K"orper. Alle auftretenden $k$-Vektorr"aume seien hier als
endlich-dimensional vorausgesetzt.
\begin{Definition}{\rm a) Eine $\mathbb R$-Filtration auf einem $k$-Vektorraum V ist
  eine monoton absteigende Abbildung
\Abb{\F^\bullet:}{\mathbb R}{\{ U; \mbox{U ist $k$-Unterraum von V} \},}{x}{\F^x}
so da"s folgende Eigenschaften erf"ullt sind:
\begin{rnumer}
\item $\F^x =V $ f"ur $x <<0$ bzw. $\F^x=(0)$ f"ur $x>>0.$
\item Sei $\F^{x-}:= \displaystyle\bigcap\limits_{y<x}\F^y.$ Dann ist $\F^{x-}=\F^{x}\; \forall x
\in \mathbb R.$
\end{rnumer}

\noindent b) Ein filtrierter Vektorraum "uber $k$ ist ein Paar $(V,\F^\bullet)$
bestehend aus einem $k$-Vektorraum V und einer $\mathbb R$-Filtration $\F^\bullet$
auf V.}
\end{Definition}

\begin{Notation}{\rm  a) In der Regel geben wir einen filtrierten Vektorraum
  $(V,\F^\bullet)$ "uber $k$ nur durch den zugrundeliegenden $k$-Vektorraum V
  an. In diesem Fall schreiben wir dann $V^x$ f"ur den Unterraum $\F^x$ von
  V.

\noindent b) Analog zum Teil a),ii) der obigen Definition setzen wir
$$ \F^{x+}:= \bigcup_{y>x} \F^y ,\; x\in \mathbb R.$$ Schlie"slich definieren wir
$$ gr^x(V):= gr^x_{\F}(V):= \F^x/\F^{x+}, \; x\in \mathbb R.$$ }
\end{Notation}

Es sei bemerkt, da"s sich alle "ublichen Operationen aus der linearen
Algebra auf filtrierte Vektorr"aume "ubertragen lassen. Die f"ur diese
Arbeit relevanten Operationen behandeln wir im nachfolgenden Beispiel.
\begin{Beispiel}{\rm  Seien  $V,W$ zwei filtrierte Vektorr"aume, so definieren wir
$$ (V \oplus W)^x := V^x \oplus W^x$$
bzw.
$$ (V \otimes W)^x:= \sum_{y+z=x} V^y \otimes W^z, \; x\in \mathbb R .$$ 
Dadurch erhalten wir sowohl auf der direkten Summe $V \oplus W$
als auch auf dem Tensorprodukt $V \otimes W$ Filtrationen.
Bezeichnet andererseits $U \subset V$ einen $k$-Unterraum von $V$, so werden durch
$$ U^x:= U \cap V^x $$
bzw.
$$ (V/U)^x:= U+V^x/U, \; x\in \mathbb R$$
auf dem Unterraum $U$ bzw. auf dem Quotienten $V/U$ Filtrationen definiert. Dabei nennen wir die letztere aus naheliegenden
Gr"unden die Quotientenfiltration.
Das sukzessive Anwenden dieser Regeln liefert ebenso Filtrationen auf 
$V^{\otimes n}, Sym^n(V)$ und $ \bigwedge^n V, n\in \mathbb N.$}
\end{Beispiel}

\begin{Bemerkung}{\rm  Ist $U$ ein Unterraum eines filtrierten Vektorraumes
  $V=(V,\F^\bullet)$,
so versehen wir $U$, sofern nichts anderes explizit gesagt wird, immer mit
der oben definierten induzierten Filtration.} 
\end{Bemerkung}

Um von der Kategorie der filtrierten Vektorr"aume "uber $k$ sprechen zu
k"onnen, m"ussen wir noch die Morphismen zwischen den jeweiligen Objekten angeben.

\begin{Definition}{\rm Ein Morphismus
$ f: (V,\F^\bullet)  \longrightarrow (V',\F'^\bullet)$
von filtrierten Vektorr"aumen "uber k, ist ein k-Homomorphismus $f: V
\rightarrow V'$ der zugrundeliegenden Vektorr"aume, so da"s $f(\F^x)
\subset \F'^x \;\forall x \in \mathbb R.$}
\end{Definition}
Wir erhalten somit die Kategorie der filtrierten Vektorr"aume "uber $k.$
Diese ist offenbar eine $k$-lineare Kategorie, welche aber nicht abelsch ist.
Trotzdem exis\-tieren kurze exakten Sequenzen in ihr. Dies ist durch die nachstehende
Definition gew"ahrleistet.
\begin{Definition}{\rm  Eine kurze exakte Sequenz 
$$ \exseq{V'}{V}{V''}$$ von filtrierten Vektorr"aumen "uber k ist eine
Folge von Morphismen, so da"s
$$ \exseq{V'^x}{V^x}{V''^x} $$ 
eine exakte Folge von k-Vektorr"aumen ist $ \forall x \in \mathbb R.$ }
\end{Definition}

Wir kommen nun zu den numerischen Invarianten filtrierter
Vektorr"aume. Sei zun"achst $V=(V,\F^\bullet)$ ein eindimensionales Objekt. Dann
existiert eine eindeutig bestimmte reelle Zahl $x_0$ mit
$gr^{x_0}(V) \neq (0).$ Diese Zahl, welche wir den Grad von $V$
nennen, bezeichnen wir mit $$\deg_{\F^\bullet}V =\deg V.$$ Allgemein hat man f"ur beliebige
filtrierte Vektorr"aume, in Analogie zu den Vektorb"undeln auf einer
Riemannschen Fl"ache,  folgende
Definition (vgl. \cite{HN} 1.1).
\begin{Definition}{\rm Sei $V=(V,\F^\bullet)$ ein filtrierter Vektorraum "uber $k$. Dann hei"st
\begin{rnumer}
\item $rk\; V:= \dim_k(V)$ der Rang von V,
\item $\deg V :=\deg_{\F^\bullet}(V):=\deg(\maxwedge{V}) $ der Grad von
  V,
\item $\slo(V):=\slo_{\F^\bullet}(V):=\frac{\textstyle \deg(V)}{\textstyle rk V}$ f"ur $V\neq (0)$ der Anstieg oder slope
  von V.
\end{rnumer} Dabei bezeichne $\maxwedge V=\bigwedge^{\dim V} V$ die maximale
"au"sere Potenz von $V$.}
\end{Definition}
\noindent Es sollen zun"achst einige einfache Eigenschaften dieser Invarianten
diskutiert werden.
\begin{Lemma} Sei $ f: V \longrightarrow W $ ein Morphismus von filtrierten
  Vektorr"aumen, der auf den zugrundeliegenden Vektorr"aumen ein
  Isomorphismus ist. Dann gilt $$\deg(V) \leq  \deg(W).$$ 
\end{Lemma}

\beweis Man reduziert den Fall sofort auf Objekte vom Rang 1, f"ur die 
die Aussage des Lemmas offensichtlich ist. \qed

\noindent Genauso wie bei den Vektorb"undeln, verh"alt sich der Grad auch in unserer
Kategorie additiv.
\begin{Lemma}  Sei \quad \exseq{V'}{V}{V''} eine kurze exakte Sequenz von
  filtrierten Vektorr"aumen. Dann
gilt $$\deg(V) = \deg(V') + \deg(V'').$$
\end{Lemma}

\beweis  Die Wahl eines Schnittes von $V \longrightarrow V''$ liefert eine 
Isomorphie ${V\cong V' \oplus V''}$ von filtrierten Vektorr"aumen. Dabei kann ein
solcher Schnitt sukzessive "uber die Komponenten der Filtrationen
von $V$ bzw. $V''$ konstruiert werden. Es folgt  $\maxwedge{V} \cong 
{\maxwedge{V'} \otimes \maxwedge{V''}}$. Offenbar folgt dann wegen
rk(\maxwedge{V'}) = rk(\maxwedge{V''}) = 1 die Behauptung unmittelbar aus
der Definition des Tensorproduktes von filtrierten Vektorr"aumen.
\qed

\begin{Korollar} Seien $U_1, U_2$ Unterr"aume eines filtrierten
  Vektorraumes V. Dann gilt
$$\deg(U_1+U_2) \geq \deg(U_1) + \deg(U_2)- \deg(U_1 \cap U_2).$$
\end{Korollar}

\beweis Betrachte die exakte Sequenz 
$$0 \longrightarrow U_1 \cap U_2 \longrightarrow U_1 \oplus U_2
\longrightarrow (U_1 + U_2)' \longrightarrow 0 $$
von filtrierten Vektorr"aumen. Dabei bezeichnet $(U_1 + U_2)'$ den
Unterraum $U_1 + U_2$ zusammen mit der Quotientenfiltration, die durch den Homomorphismus 
$U_1 \oplus U_2 \longrightarrow U_1+U_2 $ gegeben wird (vgl. Beispiel 1.3). Nun induziert
aber die Identit"at auf $U_1 + U_2$ einen Morphismus
$(U_1 + U_2)' \longrightarrow U_1 + U_2. $
Nach Lemma 1.8 gilt also $\deg((U_1 + U_2)') \leq \deg(U_1 + U_2).$ Aus
Lemma 1.9 resultiert dann aber sofort die Behauptung.
\qed

Wie schon im Fall eines eindimensionalen Objektes, beschreiben die
Sprungstellen den Grad  einer Filtration. Genauer gesagt gilt folgendes Resultat.

\begin{Proposition}
Sei V ein filtrierter Vektorraum "uber k. Dann gilt
$$ \deg V = \sum_x x\dim gr^x(V). $$
\end{Proposition}
\beweis Da sich beide Ausdr"ucke additiv bzgl. exakter Sequenzen verhalten,
ergibt sich die Behauptung leicht durch vollst"andige Induktion.
\qed

Zum Schlu"s dieses Paragraphen wollen wir noch die Definition einer
semistabilen Filtration angeben. Diese Definition ist wiederum v"ollig
analog zu dem Begriff der Semistabilit"at eines Vektorb"undels auf
einer Riemannschen Fl"ache (\cite{HN}).
Bezeichnet $W$ einen beliebigen $k$-Vektorraum und ist $K/k$ eine
K"orpererweiterung, so setzen wir zur Abk"urzung $W_K:=W\otimes_k K.$

\begin{Definition}{\rm Sei $K/k$ eine K"orpererweiterung und V ein  k-Vektorraum. 
Sei $\F^\bullet $ eine $\mathbb R$-Filtration auf $V_K.$ Die Filtration $\F^\bullet$ hei"st semistabil genau dann
wenn die Ungleichung $$\slo(U_K)  \leq \slo(V_K)$$ f"ur jeden k-rationalen
  Unterraum $U \neq (0)$ von V erf"ullt ist.}
\end{Definition}
Der Grund, da"s man in der obigen Definition K"orpererweiterungen $K/k$ betrachtet, liegt  in der Trivialit"at des Semistabilit"atsbegriffes im Fall $K=k.$
In dieser Situation sind n"amlich genau diejenigen Filtrationen semistabil, welche
genau eine Sprungstelle aufweisen.

\section{Die Untervariet"aten $\Fl_g(\B')$ und $Y(\B')$}

Von nun an setzen wir voraus, da"s unser Grundk"orper ein
endlicher K"orper ist, d.h. es ist $k=\mathbb F_q$ f"ur eine Primzahlpotenz
$q\in \mathbb N$. 

Wir fixieren einen $k$-Vektorraum $V$ der Dimension $d\geq 1$.  Au"serdem
fixieren wir eine Funktion
$$g: \mathbb R \longrightarrow \mathbb Z_{\geq 0}$$ 
mit 
$$\sum_xg(x)=d \;\mbox{ und }\; \sum_x xg(x)=0.$$ 

\noindent Aus der ersten Gleichung folgt, da"s die Funktion $g$ einen endlichen Tr"ager besitzt. 
Dem Paar $(V,g)$ l"a"st sich eine "uber $k$ definierte projektive Variet"at
$\Fl_g = \Fl_g(V)$ zuordnen (vgl. \cite{R1}, \S 2), so da"s sich f"ur jede K"orpererweiterung
$K/k$ die K-wertigen Punkte von $\Fl_g$ mit der Menge
$$\Fl_g(K)=\{ \mbox{ Filtrationen $\F^\bullet$  auf $V_K$ mit 
$\dim_K gr^x_{\F^\bullet}(V_K)=g(x)\; \forall x\in \mathbb R$}\} $$ identifiziert. 
Ebenso lassen sich die im Paragraphen 1 eingef"uhrten semistabilen
Filtrationen parametrisieren. Genauer gesagt existiert eine in $\Fl_g$
offene Untervariet"at $\Fl_g^{ss}$, so da"s f"ur jede K"orpererweiterung
$K/k$ 
$$\Fl_g^{ss}(K)=\{\F^\bullet \in \Fl_g(K) ; \F^\bullet \mbox{ ist semistabil }
   \} $$ (\cite{R1}, Prop. 2.1.).
Man nennt $\Fl_g^{ss}$ auch den zum Paar $(V,g)$ zugeh"origen {\it Periodenbereich.}
Um die Existenz dieser Periodenbereiche als Objekte in der Kategorie der Variet"aten zu sichern, hat man, wie es bei uns
der Fall ist, endliche Grundk"orper $k$ zu betrachten. F"ur solche K"orper
existieren n"amlich nur endlich viele $k$-rationale Unterr"aume $U$ in $V,$
f"ur die die Ungleichung $\slo(U_K) \leq \slo(V_K),$
bzgl. einer Filtration auf $V_K$, zu "uberpr"ufen ist.

Es sollen in dieser Arbeit au"ser $\Fl^{ss}_g$ auch noch andere
Untervariet"aten von $\Fl_g$ betrachten werden. Hierf"ur ist der Begriff einer Unterfunktion von $g$ relevant.
\begin{Definition}{\rm 
Eine Unterfunktion $h$ von $g$ ist eine Funktion \mbox{ $h : \mathbb R
  \longrightarrow \mathbb Z_{\geq 0}$ } 
mit $h(x) \leq g(x) \; \forall x \in \mathbb R.$}
\end{Definition}
Wir bezeichnen im folgenden mit  $\B=\B_g$ die Menge aller nichttrivialen
Unterfunktionen von $g,$ d.h. es gilt $\{0,g\} \cap \B = \emptyset$.
F"ur eine beliebige Funktion \mbox{$h : \mathbb R \longrightarrow \mathbb Z_{\geq 0}$} mit endlichem Tr"ager sei
$$\|h\|:= \sum\limits_x h(x)$$ als die {\it L"ange} von $h$ bezeichnet. 
Insbesondere hat die fixierte  Funktion $g$ die L"ange
$d$. F"ur eine Zahl $1\leq i \leq d-1$ sei
$$\B^i:=\{ h \in \B ; \|h\|=i \}$$
die Teilmenge der Unterfunktionen der L"ange $i$ von $g$. 
H"aufig werden wir ein Element $h \in \B^i $ durch ein ungeordnetes $i$-Tupel
$$\usupp{h}=(y_1,\ldots,y_i)$$ reeller Zahlen darstellen, welches den Tr"ager von $h$ mit
eventuellen Multiplizit"aten wiederspiegelt. 
Den gew"ohnlichen Tr"ager einer Funktion $h$ notieren wir hingegen mit
$\supp{h}.$ Gilt ferner $y_1 \geq \ldots \geq y_i ,$
so notieren wir $\usupp{h}$ auch in der Form
$$(y_1 \geq \ldots \geq y_i).$$

Als n"achstes f"uhren wir eine partielle Ordnung auf der Menge $\B$ ein.
Wir werden im n"achsten Paragraphen einen Zusammenhang dieser Ordnung zu der
Bruhatordnung einer Coxetergruppe vom Typ $A_{d-1}$ herstellen.
\begin{Definition}{\rm 
Seien $h,h'$ zwei Unterfunktionen von $g$. Setze
$$ h \geq h'  : \Longleftrightarrow \|h\| = \|h'\|$$ und aus
$$ \usupp{h}=(x_1\geq x_2\geq\ldots\geq x_{\|h\|}),\;
\usupp{h'}=(x_1'\geq x_2'\geq\ldots\geq x'_{\|h'\|}) $$ folgt $x_i \geq x_i',\; \forall i=1,\ldots,\|h\| = \|h'\|$.}
\end{Definition}

Wir ordnen nun jeder Unterfunktion $h \in \B$ bzw. jeder Teilmenge
$\B' \subset \B$ eine konstruierbare Untermenge von $\Fl_g$ zu. Hierzu
identifizieren wir im weiteren s"amtliche auftretenen Variet"aten mit ihren
abgeschlossenen Punkten. F"ur ein $h \in \B$
setzen wir
$$ Y(h):=\{ \F^\bullet \in \Fl_g ; \exists \mbox{ $k$-rationaler Unterraum $U$
von $V$ mit } \dim\, U = \|h\|$$ $$ \mbox{ und } \dim\, gr_{\F^\bullet}^x (U)=h(x)
\;\forall x\in \mathbb R \}.$$
Bezeichnet $\B' \su \B$ eine beliebige Untermenge, so seien ferner
$$Y(\B'): = \bigcup_{h \in \B'} Y(h) \;\; \mbox{ und }\;\; \Fl_g(\B'):= \Fl_g\setminus Y(\B').$$
F"ur die Mengen $Y(h)\; , h\in \B $, hat man in Analogie zu der
Abschlu"seigenschaft von Schubertzellen (vgl. \cite{J} 13.8 ) folgendes Resultat.

\begin{Proposition} 
Sei $h \in \B$. Dann gilt
$$ \overline{Y(h)} = \bigcup_{h' \geq h} Y(h').$$
\end{Proposition}
Der Beweis dieser Proposition wird im Paragraphen 3 nachgereicht. 
Zun"achst ziehen wir Konsequenzen.
\begin{Definition}{\rm Eine Teilmenge  $\B'$ von $\B$ hei"st abgeschlossen, 
falls neben  $h \in \B'$ auch alle $h' \in \B$ mit
$$ \sum_x xh'(x) \geq \sum_x xh(x)$$ in $\B'$ enthalten sind.  }
\end{Definition}
Offenbar gilt f"ur eine abgeschlossene Menge $\B'$, da"s
neben $h \in \B'$ auch alle $h' \geq h$ in $\B'$ liegen. Daher folgt aus Proposition 2.3 das folgende Korollar.
\begin{Korollar} 
F"ur eine abgeschlossene Teilmenge $\B' \su \B \,$ ist
$\, Y(\B'^i)$ abgeschlossen in $\Fl_g$, $ i=1,\ldots,d-1.$ Also ist auch
$Y(\B')= \bigcup_{i=1}^{d-1}Y(\B'^i)$ abgeschlossen in $\Fl_g$ und
$\Fl_g(\B')$ offen.
\end{Korollar}

Aus dem folgenden Lemma resultiert, da"s der Periodenbereich tats"achlich
eine der offenen Mengen $\Fl_g(\B')$ ist.
\begin{Lemma} Sei 
$$\B^{ss} =\{h \in \B ; \sum_x xh(x)>0\}.$$
Dann ist $\B^{ss}$ abgeschlossen und $\Fl_g(\B^{ss})=\Fl_g^{ss}.$
\end{Lemma}

\beweis Die Abgeschlossenheit von $\B^{ss}$ folgt unmittelbar aus der
Definition. Der abgeschlossenene Unterraum $Y(\B^{ss})$ von $\Fl_g$ l"a"st sich wie folgt beschreiben:
$$ Y(\B^{ss})=\{\F^\bullet \in \Fl_g; \exists \mbox{ $k$-rationaler Unterraum $U$ von $V$ mit } deg_{\F^\bullet}(U)>0\}.$$
Folglich ist $Y(\B^{ss})$ gerade das abgeschlossene Komplement von $\Fl_g^{ss}$ in
$\Fl_g.$ Somit gilt also $\Fl_g(\B^{ss})=\Fl_g^{ss}. $
\qed

\section{Der Beweis der Abschlu"seigenschaft}

In diesem Paragraphen beweisen wir die Proposition 2.3. Da man
diesen Beweis auf die Abschlu"seigenschaft von Bruhatzellen zur"uckf"uhrt,
m"ussen wir zun"achst einige Notationen aus der Theorie der linearen
algebraischen Gruppen einf"uhren.

Sei $T \subset GL(V)=:G$ der Diagonaltorus und $B \subset G$ die
algebraische Untergruppe der oberen Dreiecksmatrizen bzgl. einer fixierten
Basis $\{e_1,\ldots,e_d\}$ von $V$. Da die Basis offenbar $k$-rational ist,
gilt dies auch f"ur den Torus $T$ und die Borelgruppe $B$. Sei 
$$\Delta=\{\alpha_1,\ldots,\alpha_{d-1}\}$$
die Standard-Wurzelbasis zu $(T,B)$, d.h. es gilt
$$\alpha_i (\mbox{diag}(t_1.\ldots,t_d) )=t_it_{i+1}^{-1},\; 1\leq i\leq d-1,$$wobei $\mbox{diag}(t_1,\ldots,t_d)$ die Diagonalmatrix mit
Eintr"agen $t_1,\ldots,t_d$ sei.
Wir bezeichnen mit  $W \cong S_d$ die Weylgruppe von $G$ bzgl. $T$ und mit
$$S=\{s_1,\ldots, s_{d-1}\}$$ die Menge der einfachen Spiegelungen, wobei
$s_i$ die zur Wurzel $\alpha_i$ korrespondierende Spie\-gelung ist.
F"ur eine Teilmenge $I\subset S$ sei 
$$W_I=<s_i  \; ; s_i \in I>$$ die von $I$ erzeugte
parabolische Untergruppe in $W$. Entsprechend notieren wir mit 
$$P_I=BW_IB$$ die zugeh"orige
std.-parabolische Untergruppe von $G.$ Als Extremf"alle hat man also
$P_\emptyset=B$ bzw. $P_S=G$. Ist speziell $I$ von der Gestalt
$I=S\setminus\{s_i\}$ f"ur ein $i$, so schreiben wir auch $W_{\hat{i}}$ an
Stelle von $W_I$ bzw.  $P_{\hat{i}}$ an Stelle von $P_I.$  

\medskip
Zu unserer fixierten Funktion $g$ mit
$$ \usupp{g}=(x_1 \geq \ldots \geq x_d)  $$ assoziieren
wir den reellen Cocharakter
$$\mu=\mu_g=(x_1,\ldots,x_d) \in X_{\ast}(T)_{\mathbb
  R}:=X_{\ast}(T)\otimes_{\mathbb Z} \mathbb R. $$ 
Wir bezeichnen mit 
$$W_\mu=\mbox{Stab}_W(\mu)$$ den Stabilisator von $\mu$ bzgl. der Operation
von $W$ auf $X_{\ast}(T)_{\mathbb R}$. Analog sei 
$$P_\mu=BW_\mu B$$ die
zugeh"orige parabolische Untergruppe. 
Schreibt sich der herk"ommliche Tr"ager von $g$ in der Form  supp$(g)=\{x_1',\ldots, x_r'\},$ so besteht $P_\mu$ bzgl. unserer
fixierten Basis also aus Matrizen der Form

$$\left(\matrix{M_1&*&\cdot&\cdot&*\cr
                0&M_2&\cdot&\cdot&*\cr
                \cdot&0&\cdot&\cdot&\cdot\cr
                \cdot&\cdot&\cdot&\cdot&*\cr
                0&0&\cdot&0&M_r} \right) , $$

\noindent wobei die $M_i$ invertierbare $g(x_i')\times g(x_i')$-Matrizen darstellen.
Entsprechend ist die korrespondierende parabolische Untergruppe $W_\mu$ von
$W$ isomorph zu
$$ S_{g(x_1')}\times S_{g(x_2')}\times\cdots \times S_{g(x_r')}.$$

Wir wollen nun die in Paragraph 2 eingef"uhrte Halbordnung auf $\B$ analysieren. 
F"ur ein Weylgruppenelement $w \in W$ und ein $0 \leq i \leq d$
definieren wir die folgende Unterfunktion der L"ange $i$ von $g$:
$$h_w^i(x)=\#\{k; x=(w\mu)_k, k=1,\ldots,i\}, \; x \in \mathbb R.$$
Dabei verstehen wir unter $(w\mu)_k$ den $k$-ten Eintrag im Tupel $w\mu.$
Folglich gilt dann
\begin{equation}\usupp{h_w^i}=(x_{w^{-1}(1)},\ldots,x_{w^{-1}(i)}).
\end{equation}
Als Extremf"alle hat man
$$ h^0_w=0 \;\; \forall w \in W \mbox{ bzw. } h^d_w=g \;\; \forall w \in W.$$
Wir erhalten f"ur jedes $1 \leq i \leq d-1$ eine Abbildung
$$\tilde{\kappa_i}: W \longrightarrow \B^i. $$ 
\begin{Lemma}
Die Abbildung $\tilde{\kappa_i}$ ist surjektiv und faktorisiert  "uber
$W_{\hat{i}}\backslash W/W_\mu.$ Die induzierte Abbildung
$$\kappa_i:W_{\hat{i}}\backslash W/W_\mu \longrightarrow \B^i $$ ist eine Bijektion, welche ordnungsumkehrend ist, wenn $W_{\hat{i}}\backslash W/W_\mu$
mit der Bruhatordnung versehen ist. 
\end{Lemma} 

\beweis Die erste Aussage des Lemmas ist trivial. Die
Injektivit"at der Abbildung $\kappa_i$ ist ebenso leicht einsehbar. 
Seien nun $w,w'$ zwei Kostant-Repr"asentanten bzgl. 
$W_{\hat{i}}\backslash W.$ Es gilt  
(vgl. \cite{FZ} 3.2.)
$$ w \leq w' \,\aequi\, w^{-1} \leq w'^{-1} \,\aequi\, w^{-1}(j) \leq w'^{-1}(j) \;\,
j=1,\ldots, i.$$ Hieraus ergibt sich zusammen mit der Ungleichungskette
$x_1 \geq \ldots \geq x_d$ die Behauptung.
\qed

\begin{Notation}{\rm  Wir bezeichnen f"ur ein $h \in \B^i, 1\leq i \leq d-1,$ 
mit $$ w_h \in W_{\hat{i}}\backslash W/W_\mu$$ 
das Urbild von $h$  unter $\kappa_i.$
Den Kostant-Repr"asentanten der Doppelnebenklasse $w_h$ notieren wir
mit $\dot{w_h}.$}
\end{Notation}

F"ur einen $k$-rationalen Unterraum $U$ der Dimension $i$ von $V$ und ein Element $w \in  W_{\hat{i}}\backslash W/W_\mu$ betrachten wir die zugeh"orige verallgemeinerte "`Schubertzelle"'
$$ Y_U(w)=\{\F^\bullet \in \Fl_g ;
\dim gr_{\F^\bullet}^x(U)=h^i_w(x) \; \forall x \in \mathbb R \}.$$ 
Wir wollen im nachfolgenden eine andere Beschreibung f"ur $Y_U(w)$
herleiten. Es wird sich herausstellen, da"s diese Menge einen homogenen Raum
unter einer parabolischen Untergruppe von $G$ darstellt. Betrachte hierzu die vollst"andige $k$-rationale Standardflagge
$$ V^{\bullet} = (0) \su V^1 \su V^2 \su \cdots \su V^{d-1} \su V $$
mit 
$$V^i:=<e_1,\ldots,e_i>,\; i=1,\ldots,d-1.$$
Sei  $\tilde{g}: \mathbb R \longrightarrow \mathbb Z_{\geq 0}$ definiert durch 
$\tilde{g}(x):=\sum_{y \geq x }g(y), x\in \mathbb R. $ Wir definieren einen
$k$-rationalen Punkt $\tilde{V}^{\bullet} \in \Fl_g$ durch
$$\tilde{V}^x:= <e_1,\ldots,e_{\tilde{g}(x)}>, \,x \in \mathbb RR.$$
Die nat"urliche Operation von  $P_{\hat{i}}$ auf $\Fl_g$ l"a"st  $Y_{V^i}(w)$
invariant. Wir werden in der nachstehenden Proposition sehen, da"s
$Y_{V^i}(w)$ bereits einen $P_{\hat{i}}$-Orbit in der Flaggenvariet"at
$\Fl_g$ darstellt.

\begin{Proposition} Sei $w \in W_{\hat{i}}\backslash W/W_\mu$. Dann operiert  
$P_{\hat{i}}$ transitiv auf $Y_{V^i}(w).$
Stellen wir unsere Flaggenvariet"at $\Fl_g$ als homogenen Raum $G/P_\mu$ mit
Basispunkt $\tilde{V}^{\bullet}$ dar, so gilt $Y_{V^i}(w)=P_{\hat{i}}wP_\mu /P_\mu.$ 
\end{Proposition}

\beweis Zun"achst notieren wir die zu $g$ assoziierte 1-PS $\mu$ wieder in der Form
$$\mu =(x_1,\ldots,x_d).$$ Sei $w' \in W$ ein beliebiges Element. Dann
gilt f"ur jedes $1 \leq j \leq d-1$
\[ h^j_{w'} (x)-h^{j-1}_{w'} (x) = \left\{ 
\begin{array}{r@{\quad: \quad}l} 1 & x=x_{w'^{-1}(j)} \\ 0 &
\mbox{sonst} \end{array} \right. .\]

\noindent Wir definieren f"ur jedes $w' \in W$ die Menge
\begin{eqnarray*}Y(w') & := & \{ \F^\bullet \in \Fl_g; \dim gr^x_{\F^\bullet} (V^j / V^{j-1})
  =  h^j_{w'}(x)-h^{j-1}_{w'} (x) \; \\ \\ & & \;\;\;\;\; \forall x \in \mathbb R,  j=1,\ldots,d-1 \} \\ \\
&  = & \{\F^\bullet \in \Fl_g; gr^{x_{w'^{-1}(j)}}_{\F^\bullet}(V^j/V^{j-1}) \neq (0), \;
j=1,\ldots d-1 \}.
\end{eqnarray*}
\noindent Wir ben"otigen die folgenden Lemmata.
\begin{Lemma} Unter der Identifikation $\Fl_g=G/P_\mu$ gilt
$$Y(w') = Bw'P_\mu / P_\mu.$$
\end{Lemma}

\beweis In der Tat, die Bruhatzelle $B w' P_\mu /P_\mu$ l"a"st sich  unter
der obigen Identifikation wie folgt beschreiben (\cite{LM} 1.2):
$$B w' P_\mu / P_\mu=\{\F^\bullet \in \Fl_g;  gr^{x_j}_{\F^\bullet}(V^{w'(j)}/V^{w'(j)-1})
\neq 0 \;\;  j=1,\ldots,d\}.$$ 
Die Behauptung ergibt sich nun durch eine Variablentransformation.\qed

\begin{Lemma} Die $Y(w')$ induzieren eine Zerlegung
$$Y_{V^i}(w)=\bigcup_{w' \in W_{\hat{i}}w} Y(w').$$
\end{Lemma}

\beweis Sei $\F^\bullet \in Y(w')$ f"ur ein  $w' \in W_{\hat{i}}w.$ Dann gilt
\begin{eqnarray*}\dim gr^x_{\F^\bullet} (V^i) &=& \sum_{j=1}^i \dim gr^x_{\F^\bullet} (V^j
  /V^{j-1}) = \sum_{j=1}^i \left(h^j_{w'} (x)-h^{j-1}_{w'} (x)\right) \\ & =&  h^i_{w'}(x) =  h^i_w (x).
\end{eqnarray*}
Sei nun umgekehrt  $\F^\bullet$ ein Element aus $Y_{V^i}(w)$. Da s"amtliche
Subquotienten $V^j/V^{j-1}, j=1,\ldots,d-1$ eindimensional sind,
existiert zu jedem $1 \leq j \leq d-1$ genau eine Sprungstelle
$y_j \in \supp{g} \mbox{ mit } gr_{\F^\bullet}^{y_j}(V^j/V^{j-1})\neq (0).$
Nun gilt aber
$\dim gr_{\F^\bullet}^x(V^i)=h^i_w(x)\; \forall x\in \mathbb R . $ Folglich sind die
ersten $i$ Sprungstellen $y_1,\ldots,y_i$ lediglich eine Permutation der 
$x_{w^{-1}(1)}, \ldots, x_{w^{-1}(i)}.$ Also existiert ein $w' \in
W_{\hat{i}}w$ mit $\F^\bullet \in Y(w').$ 
\qed
 
{\noindent \bf Ende des Beweises von Proposition 3.3:} Die
Identit"at $P_{\hat{i}}=BW_{\hat{i}}B=\bigcup_{w\in W_{\hat{i}}}BwB$
  induziert die Zerlegung (\cite{B} Kap. IV, \S 2 Prop. 2)
$$P_{\hat{i}}w P_\mu /P_\mu = \bigcup_{w' \in  W_{\hat{i}}w} Bw' P_\mu /
P_\mu .$$ 
Die Behauptung der Proposition 3.3 ergibt sich nun unmittelbar aus den
vorhergehenden Lemmata. \qed

\begin{Bemerkung}{\rm  Die Bezeichnung verallgemeinerte "`Schubertzelle"' ist
eigentlich irref"uhrend, da die Variet"aten $P_{\hat{i}}wP_\mu /P_\mu,$
im Gegensatz zu den her\-k"ommlichen Bruhatzellen, im allgemeinen keine standard-affinen R"aume sind.}
\end{Bemerkung}  

\medskip
Man hat also eine Zerlegung
$$ \Fl_g = \stackrel{\cdot}{\bigcup\limits_{w \in W_{\hat{i}}\backslash W/W_\mu}}Y_{V^i}(w) 
= \stackrel{\cdot}{\bigcup\limits_{w \in W_{\hat{i}}\backslash W/W_\mu}}P_{\hat{i}}\ w P_\mu/P_\mu  $$
in verallgemeinerte "`Schubertzellen"', die der Zerlegung in
$P_{\hat{i}}$-Orbiten von $G/P_\mu$,
entspricht. Genauso induziert jeder $k$-rationale Unterraum $U$ der
Dimension $i$ von $V$ eine
disjunkte Vereinigung
\begin{equation}\Fl_g = \stackrel{\cdot}{\bigcup\limits_{w \in
      W_{\hat{i}}\backslash W/W_\mu}}Y_U(w).
\end{equation}
Genauso wie im
Fall von herk"ommlichen Bruhatzellen haben wir auch bei unseren 
verallgemeinerten "`Schubertzellen"' eine entsprechende
Abschlu"seigenschaft. 

\begin{Lemma}
F"ur ein beliebiges $w\in W_{\hat{i}}\backslash W
/W_\mu, 1\leq i \leq d-1$ gilt 
$$\overline{ P_{\hat{i}}wP_\mu /P_\mu}= \bigcup_{w' \in
  W_{\hat{i}}\backslash W /W_\mu \atop w' \leq w} P_{\hat{i}}w'P_\mu /P_\mu$$ 
bzw.
$$\overline{ Y_{V^i}(w)}= \bigcup_{w \in W_{\hat{i}}\backslash W /W_\mu
  \atop w' \leq w} Y_{V^i}(w').$$
\end{Lemma}

\beweis Der Beweis folgt aus der Abschlu"seigenschaft f"ur herk"ommliche
Bruhatzellen  und aus der Gleichheit
$$ P_{\hat{i}}wP_\mu = BW_{\hat{i}}\dot{w}P_\mu = \bigcup_{w' \in
  W_{\hat{i}}\dot{w}} Bw'P_\mu  $$ 
(vgl.  \cite{B} Ch. IV, \S 2 Prop. 2).
\qed

Wir wollen  nun  die Proposition 2.3 beweisen. Zuvor sei aber noch f"ur
$1\leq i \leq d-1$ die Identit"at 
\begin{equation}
Y(h^i_w) = \bigcup_{g \in G(k)} g Y_{V^i}(w)
\end{equation}
erw"ahnt. Sie folgt aus der Tatsache, da"s die endliche Gruppe
$G(k)$ transitiv auf den $k$-rationalen Unterr"aumen einer gegebenen Dimension operiert.

\bigskip
{\noindent \bf Beweis von Proposition 2.3:} Sei $h \in \B^i$. 
Dann gilt wegen (3) die Gleichheit
$$ Y(h)=\bigcup_{g \in G(k)}g Y_{V^i}(w_h).$$
Aus der Abschlu"seigenschaft der verallgemeinerten "`Schubertzellen"'
folgt schlie"slich
\begin{eqnarray*} \overline{Y(h)} &=& \bigcup_{g \in
    G(k)}\overline{gY_{V^i}(w_h)}= \bigcup_{g\in G(k)} \bigcup_{w \in
    W_{\hat{i}}\backslash W /W_\mu \atop  w \leq w_h} gY_{V^i}(w) \\
 &=& \bigcup_{w \in W_{\hat{i}}\backslash W /W_\mu \atop  w \leq w_h}
    Y(h^i_{w'})  = \bigcup_{h'\geq h} Y(h') \;\;\;\;\square.
\end{eqnarray*}

\section{Die Formulierung der Hauptresultate}

Sei $\ell$ eine f"ur den Rest der Arbeit fixierte Primzahl mit $(\ell,q)=1.$
Bezeichnet $X$ eine Variet"at "uber dem endlichen K"orper $k$ und $n\in \mathbb N$ eine nat"urliche
Zahl, so bezeichnen wir zur Abk"urzung seine $\ell$-adische Kohomologie mit
$$ H^n_{\acute{e}t}(X):=H^n_{\acute{e}t}(X \times_k \overline{k}, \mathbb Q_\ell). $$
Dies ist ein $Gal(\overline{k}/k)$-Modul, das hei"st, er ist mit einer
Operation des Frobenius versehen.
Die $\ell$-adische Kohomologie mit kompakten Tr"ager von 
$X \times_k \overline{k}$  notieren wir mit 
$$H^n_c(X):=H^n_c(X\times_k \overline{k},\mathbb Q_\ell).$$ 
Ist schlie"slich $V$ ein beliebiger $\mathbb Q_\ell$-Vektorraum und $r$ eine
ganze Zahl, so notieren wir den $r$-fachen Tate-Twist von $V$ mit 
$V(r)=V\otimes_{\mathbb Q_\ell} \mathbb Q_\ell(r).$

Sei nun $I \subset S$ eine echte maximale Untermenge,
d.h. es ist $I=S\setminus\{s_i\}$ f"ur ein $1 \leq i \leq d-1.$ Wir definieren
\begin{equation} \Omega_I:=\bigcup\limits_{w \in W \atop h^i_w \in \B'^i}
W_{\hat{i}}w W_\mu \cap W^\mu .
\end{equation}
Einer beliebigen Teilmenge $I \subn S $ ordnen wir dann die Menge
\begin{equation} \Omega_I:= \bigcap_{I\subset J \atop \#(S\setminus J)=1}
  \Omega_J.
\end{equation}
zu. Schlie"slich setzen wir  $\Omega_S= W^\mu.$
Offenbar ist die Zuordnung $I \mapsto \Omega_I$ inklusionserhaltend, genauer gilt sogar die Gleichheit
\begin{equation} \Omega_{I \cap J} = \Omega_I \cap \Omega_J \;\;\; \forall I,J
  \subset S.
\end{equation}
Im weiteren bezeichnen wir  f"ur ein Weylgruppenelement $w \in W^\mu$ mit
$I_w$ die kleinste Teilmenge von $S,$ so da"s $w$ in $\Omega_{I_w}$ enthalten
ist. 
Offenbar gilt dann die "Aquivalenz
\begin{equation} I_w \subset I \; \aequi \;w \in \Omega_I.\end{equation}

F"ur eine parabolische Untergruppe $P\subset G$ betrachten wir die triviale
Darstellung von $P(k)$ auf $\mathbb Q_\ell.$ Die resultierende induzierte Darstellung
$i_{P(k)}^{G(k)}(\mathbb Q_\ell)$ notieren wir zur Abk"urzung mit $i_P^G.$ 
Ferner setzen wir $$v^G_P=i^G_P/\sum_{P\subn Q}i^G_Q.$$ Im Fall $P=B$
erhalten wir somit die Steinberg-Darstellung \cite{S}. Ist schlie"slich $w\in W$ ein Weylgruppenelement, so bezeichnen wir mit
$l(w)$ seine L"ange.

Mit diesen Notationen k"onnen wir nun die Hauptresultate dieser
Arbeit formulieren.
\begin{Satz} Sei $\B' \subset \B^{ss}$ eine abgeschlossene Teilmenge. Dann gilt
$$H^{\ast}_{\acute{e}t}(Y(\B'))=\bigoplus_{w\in W^\mu \atop
\#(S\setminus I_w)=1} i_{P_{I_w}}^G(-l(w))[-2l(w)]\oplus $$
$$\bigoplus_{w\in W^\mu \atop \#(S\setminus I_w)>1} 
\left(i^G_G(-l(w))[-2l(w)]\oplus
v_{P_{I_w}}^G(-l(w))[-2l(w)-\#(S\setminus I_w)+1]\right).$$
{\rm Die Notation $[-n], n\in \mathbb N,$ bedeutet dabei, da"s der
voranstehende Vektorraum in den Grad $n$ des betreffenden graduierten Kohomologieringes geshiftet wird.}
\end{Satz}
F"ur das offene Komplement $\Fl_g(\B')$ von $Y(\B')$ hat man das folgende Resultat.

\begin{Satz}
Sei $\B' \subset \B^{ss}$ abgeschlossen. Dann gilt
 $$H^{\ast}_c(\Fl_g(\B'))=\bigoplus_{w \in W^\mu}
v^G_{P_{I_w}}(-l(w))[-2l(w)-\#(S\setminus I_w)].$$
\end{Satz}

Ist speziell $\B' = \B^{ss},$ so l"a"st sich die Menge $I_w$ wie folgt
beschreiben. Bezeichne mit 
$$\langle\;\,,\;\rangle: X_\ast(T)_{\mathbb R} \times X^\ast(T)_{\mathbb R} \longrightarrow \mathbb R $$ 
die nat"urliche Paarung zwischen den reellen Cocharakteren und Charakteren.
Sei $\alpha_i^\vee \in X_\ast(T)_{\mathbb R}$ die Cowurzel zu $\alpha_i.$ 
Sei schlie"slich  $\omega_i$ das zur Wurzel $\alpha_i$ korrespondierende
Fundamentalgewicht, d.h. es gilt
$\langle \alpha_i^\vee,\omega_j \rangle=\delta_{ij}.$

\begin{Lemma}
F"ur jedes $w \in W^\mu$ gilt
$$ I_w=\{ s_i \in S; \langle w\mu, \omega_i\rangle \leq 0\} .$$
\end{Lemma}
\beweis Aus (4) und (7) folgt die Identit"at 
\begin{eqnarray*}I_w &=& \bigcap_{w \in \Omega_I \atop
  \#(S\setminus I) =1} I = S \setminus \{ s\in S ; w \in \Omega_{S\setminus
  \{s\}} \} 
= S \setminus \{s_i \in S; h^i_w \in (\B^{ss})^i\}\\
&=& \{ s_i \in S; h^i_w \not\in (\B^{ss})^i \}.
\end{eqnarray*}
Per Definition von $\B^{ss}$ gilt aber $h\in \B^{ss}\; \aequi\; \sum_xxh(x)>0.$
Aus (1) erh"alt man die Gleichung $\sum_x xh_w^i(x) =
\langle w\mu,\omega_i\rangle$ und somit die Behauptung.\qed

\noindent Definieren wir f"ur jedes $w \in W^\mu$ die Menge $\Delta_w$ durch
$$\Delta_w:=\{\alpha_i\in \Delta; s_i \not\in I_w \}, $$
so ist $P_{I_w}$ gerade diejenige std.-parabolische Untergruppe, deren
unipotentes Radikal von $\Delta_w$ erzeugt wird. Als Spezialfall des Satzes
4.2 erh"alt man somit folgende Formel, welche von Kottwitz und Rapoport
vermutet wurde, vgl. \cite{R3}.
\begin{Satz} Es gilt
 $$H^{\ast}_c(\Fl_g^{ss})=\bigoplus_{w \in W^\mu}
v^G_{P_{I_w}}(-l(w))[-2l(w)-\#\Delta_w].$$
\end{Satz}

Als n"achstes leiten wir folgenden Verschwindungssatz her (vgl. \cite{R3}).
\begin{Korollar}
F"ur den zu $g$ assoziierten Periodenbereich $\Fl^{ss}_g$ ist
$$ H^i_c(\Fl_g^{ss})=0 \;\;\; 0 \leq i \leq d-2$$
und
$$ H^{d-1}_c(\Fl_g^{ss})= v^G_B .$$
\end{Korollar}
F"ur den Beweis  ben"otigen wir nachstehende Lemmata.
\begin{Lemma}
Sei $w\in W^\mu$ und $s=s_\alpha \in S$ die zu einer Wurzel $ \alpha \in
\Delta$ zugeh"orige Spiegelung. Wir setzen
voraus, da"s $sw$ in $W^\mu$ enthalten und die Identit"at
$l(sw)=l(w)+1$ erf"ullt ist.  Dann gilt
$$ \Delta_{sw}\setminus \{\alpha \} = \Delta_w \setminus \{ \alpha\} $$ 
und
$$ \Delta_{sw} \subset \Delta_w .$$
Insbesondere unterscheiden sich $\Delta_{sw}$ und $\Delta_w$ um h"ochstens
ein Element.
\end{Lemma}

\beweis Per Definition von $\Delta_w$ gilt
$$ \Delta_w= \{\alpha \in \Delta; \langle w\mu,\omega_\alpha \rangle\;> 0\}. $$ 
Ist nun  $\beta \in \Delta$ beliebig, so gilt
$ \langle sw\mu, \omega_\beta\rangle\;=\; \langle w\mu,s\omega_\beta\rangle $
und
$$ s\omega_\beta=\left\{\matrix{\omega_\beta &;& \beta \neq \alpha \cr
\omega_\alpha -\alpha &;& \beta=\alpha } \right. \;\;(\cite{B}\; 1.10).$$
Hieraus folgt offenbar die erste Behauptung. Die Voraussetzung
$l(sw)=l(w)+1$ ist "aquivalent zur Positivit"at der Wurzel $w^{-1}(\alpha),$  d.h.
$w^{-1}(\alpha) > 0 \;\;(\cite{H}\; 1.6).$ Da $\mu$ im Abschlu"s
der positiven Weylschen Kammer liegt (vgl. \S 3 ) folgt 
$\langle w\mu,\alpha \rangle\; =\; \langle \mu,w^{-1}\alpha \rangle\; \geq
0.$ Ist  $\alpha$ ein Element von $\Delta_{sw},$ so
gilt also 
$$ 0< \langle sw\mu,\omega_\alpha\rangle\;=\;\langle w\mu,\omega_\alpha - \alpha\rangle\;=\;
\langle w\mu,\omega_\alpha\rangle  - \langle w\mu,\alpha \rangle .$$ Es folgt
$ \langle w\mu, \omega_\alpha\rangle \;  > \; \langle
w\mu,\alpha\rangle\; \geq 0,$ also  $\alpha \in \Delta_w.$ \qed

\begin{Lemma}
Sei $w\in W^\mu.$ Dann gilt
$$ \#(\Delta \setminus \Delta_w) \leq l(w).$$
\end{Lemma}
\beweis Sei $w=s_1\cdots s_r,\; r=l(w)$ eine reduzierte Darstellung mit $s_i
\in S$. Wegen
der Gleichheit 
$$ W^\mu=\{ w \in W; l(ws) = l(w)+1 \; \forall s \in W_\mu\cap S \}$$
(\cite{H} 1.10) sind neben $w$ auch s"amtliche Ausdr"ucke der Form
$$ s_r,\; s_{r-1}s_r, \ldots, s_2\cdots s_r$$
in $W^\mu$ enthalten. Die Behauptung ergibt sich nun durch vollst"andige
Induktion aus  Lemma 4.6.\qed

{\noindent \bf Beweis von Korollar 4.5:} Die Anwendung von Lemma 4.7
liefert f"ur jedes $w\in W^\mu$ die Ungleichung
$$ \# \Delta_w +2l(w)\geq (\# \Delta -l(w)) +2l(w)=\# \Delta +l(w) \geq \#
\Delta.$$ Gleichheit gilt dabei in dieser Kette genau dann wenn $w=1.$ 
Durch Anwendung von Satz 4.4 erh"alt man unter Ber"ucksichtigung von $\#\Delta = d-1$ die Behauptung.\qed

Wir wollen nochmals die Gegenl"aufigkeit der Ausdr"ucke $l(w)$ und
$\Delta_w$ bez"uglich der Bruhatordnung auf $W^\mu$ betonen. Aus $w' \leq
w$ folgt $l(w') \leq l(w),$ aber $\Delta_w \subset \Delta_{w'}.$ Die
letztere Inklusion sieht man folgenderma"sen ein. Gilt $w' \leq w,$ so hat
man nach \cite{FZ} 3.2. die Relation $w'^{-1} \omega_i \geq
w^{-1}\omega_i \;\forall i.$ Die obige Inklusion folgt nun aus der Tatsache,
da"s $\mu$ im Abschlu"s der positiven Weylschen Kammer liegt.
Daher kann es passieren, da"s f"ur $w' \leq w$ der induzierte Grad von $w'$
gr"o"ser als derjenige von $w$ ist, wie folgendes Beispiel zeigt. Dies
widerlegt eine Behauptung von Rapoport in \cite{R3}.
\begin{Beispiel} {\rm Sei $\dim_kV=5$  und $g:\mathbb R \rightarrow \mathbb
    Z_{\geq 0}$ eine Funktion mit geordneten Tr"ager
$$\usupp{g}=(x_1>x_2>x_3>x_4>x_5),$$
wobei $\sum_{i=1}^5 x_i=0$ und $x_4>0.$ Betrachte die Weylgruppenelemente
$$w'=(2,3,4) \mbox{ bzw. } w=(2,3,4,5).$$ 
Offenbar ist $w' < w,$ und es gilt $l(w')=2$ bzw. $l(w)=3.$ Weiter ist
$$ \Delta_{w'}=\{\alpha_1,\alpha_2,\alpha_3,\alpha_4\} \mbox{ bzw. }
\Delta_w=\{\alpha_1\}.$$
Somit gilt also
$$ 2l(w') + \#\Delta_{w'}= 8\; \mbox{ aber }\; 2l(w)+\#\Delta_w=7. $$}
\end{Beispiel}

Wir wollen nun die Reinheit der Kohomologie von Periodenbereichen diskutieren.
Wie man der Formel aus Satz 4.4 entnimmt, ist die Kohomologie der
Periodenbereiche im allgemeinen nicht rein. 
\begin{Beispiel}{\rm Sei $\dim_kV=3$ und $g:\mathbb R \rightarrow \mathbb
  Z_{\geq 0} $ eine Funktion mit $\usupp{g}=\{x_1>x_2>x_3\},$ wobei $x_2 >0$
und $x_1+x_2+x_3=0.$ In diesem Fall gilt f"ur  die Kohomologie von $\Fl_g^{ss}$

$\begin{array}[t]{lcl} \\
 H_c^0(\Fl^{ss}_g)=0 &\;\;\;\;\;\;\;\;\;\; & H_c^4(\Fl^{ss}_g)=v^G_B(-1)\oplus i^G_G(-2) \\ \\
 H_c^1(\Fl^{ss}_g)=0 &\;\;\;\;\;\;\;\;\;\; & H_c^5(\Fl^{ss}_g)=v^G_{P_{\{s_1\}}}(-2) \\ \\
 H_c^2(\Fl^{ss}_g)=v^G_B &\;\;\;\;\;\;\;\;\;\; & H_c^6(\Fl^{ss}_g)=i^G_G(-3) \\ \\
 H_c^3(\Fl^{ss}_g)=v^G_{P_{\{s_1\}}}(-1) &\;\;\;\;\;\;\;\;\;\; & H_c^i(\Fl^{ss}_g)=0 \;\; \forall
i>6
\end{array}$}
\end{Beispiel}
Dagegen ist die Kohomologie des Drinfeldraumes $\Omega(V)$ im
 folgenden Sinn rein.  Dabei ist $\Omega(V) \subset \mathbb P(V)$ gleich
$\Fl_g^{ss}$ f"ur $g$ mit $\usupp{g}=(x_1,x_2^{dim V -1}),\; x_1 > x_2.$ 
Explizit ist
$\Omega(V)$ das Komplement aller $k$-rationalen Hyperebenen im
$\mathbb P(V),$ d.h. es gilt 
$$ \Omega(V)=\mathbb P(V)\setminus \bigcup_{H \subset V \;\;\mbox{\tiny  $k$-rat. } \atop  \dim H=\dim V-1}\mathbb P(H).$$
In diesem Fall ist die Menge $W^\mu$ durch die Elemente
$$w_0=1,w_1=s_1,w_2=s_2s_1,\ldots,w_{d-1}=s_{d-1}s_{d-2}\cdots s_1$$
gegeben. Wegen der Gleichheit
$$\Delta_{w_i}=\{\alpha_{i+1},\ldots,\alpha_{d-1}\}.$$
ergibt sich f"ur die Kohomologie des Drinfeldraumes folgende explizite Formel:
$$H^{\ast}_c(\Omega(V),\mathbb Q_\ell)=\bigoplus_{i=0}^{d-1}
v^G_{P_{(d-i,1,\ldots,1)}}(i+1-d)[-2(d-1)+i].$$
Also  ist $H_c^{2(\dim V-1)-i}(\Omega(V))$ ein reiner
Galois-Modul vom Gewicht $2(\dim V-1)-2i,\;
i=0,\ldots,2(\dim V-1).$ Diese Reinheit ist hier eine Folgerung unseres
Hauptresultates. Sie ist allerdings a priori klar, da dies allgemein f"ur das Komplement eines
Hyperebenenarrangements gilt (\cite{OT} S. 194). Diese Tatsache
erm"oglicht die Berechnung der einzelnen Kohomologiegruppen des
Drinfeldraumes aus der Formel f"ur die Euler-Poncar\'e-Chrakteristik, da sich die Glieder verschiedenen Grades nicht gegenseitig wegk"urzen k"onnen. 

Zum Abschlu"s dieses Paragraphen wollen wir die Beweisstrategie der S"atze
4.1 bzw. 4.2 schildern. Im Fall von Satz 4.1 betrachten wir im Paragraphen 5
eine "Uberdeckung von $Y(\B')$ durch abgeschlossene Untervariet"aten der
Gestalt $gY_I$, wobei $I$ die  maximalen echten Teilmengen von $S$
und $g$ die Linksnebenklassen in $G(k)/P_I(k)$ durchl"auft. "Ahnlich wie
bei der Herleitung der Mayer-Vietoris-Sequenz bzgl. einer abgeschlossenen "Uberdeckung konstruieren wir f"ur eine \'etale Garbe $F$ auf $Y(\B')$ eine Sequenz von \'etalen Garben auf
diesem Raum. Anstatt  s"amtliche sukzessiven Durchschnitte der
Untervarit"aten $gY_I$ zu betrachten, wie im Fall der
Mayer-Vietoris-Sequenz,  nehmen wir hier nur Durchschnitte
der Form 
$$g(Y_{I_1} \cap \ldots \cap Y_{I_r}),\;\; g\in G(k)/P_{I_1 \cap \ldots
  \cap I_r}(k),$$ wobei $I_1,\ldots I_r$  maximale echte Teilmengen
in $S$ sind. Auf diese Weise erhalten wir einen Komplex 
$$0 \rightarrow F \rightarrow \bigoplus_{I \subset S \atop
\#(S\setminus I)=1} \bigoplus_{g \in (G/P_I)(k)} (\phi_{g,I})_*(\phi_{g,I})^* F \rightarrow
\bigoplus_{I \subset S \atop \#(S\setminus I)=2} \bigoplus_{g \in (G/P_I)(k)}
(\phi_{g,I})_*(\phi_{g,I})^* F \rightarrow $$
$$ \dots \rightarrow \bigoplus_{I \subset S \atop \#(S\setminus I)=d-2} \bigoplus_{g \in (G/P_I)(k)} (\phi_{g,I})_*(\phi_{g,I})^* F \rightarrow
\bigoplus_{g \in (G/B)(k)} (\phi_{g,\emptyset})_*(\phi_{g,\emptyset})^* F
\rightarrow 0$$ von \'etalen Garben auf $Y(\B'),$
 wobei die $\phi_{g,I}$ die abgeschlossenen Immersionen 
$gY_I \hookrightarrow Y$ bezeichnen.
Die  charakterisierende Indexmenge dieses Komplexes entspricht dabei dem
Tits-Komplex zur Gruppe $G(k).$  Im Paragraphen 6 zeigen wir die
Azyklizit"at dieses Komplexes, indem wir die Azyklizit"at des durch
Lokalisierung in einem geometrischen Punkt $x$ von $Y$ entstehenden
Komplex\-es beweisen. Dieser entspricht einem Kettenkomplex
mit Werten in $F_x,$ der von einem Unterkomplex des Tits-Komplexes
induziert wird. Mit Hilfe eines Kontraktionskriteriums f"ur simpliziale Komplexe,
die durch partiell geordnete Mengen induziert werden, zeigen wir die
Zusammenziehbarkeit der geometrischen Realisierung dieses Unterkomplexes. 
Hieraus folgt insbesondere die Azyklizit"at des
betrachteten Kettenkomplexes. Im Paragraphen 7 werten wir schlie"slich die
Spektralsequenz aus, die sich aus dem Ausgangskomplex ergibt. 
Dies liefert
die gesuchten Kohomologiegruppen von $Y.$
Den Satz 4.2 beweisen wir durch die Anwendung der langen exakten
Kohomologiesequenz zu dem Tripel $\Fl_g(\B') \hookrightarrow \Fl_g
\hookleftarrow Y(\B').$

\section{Der fundamentale Komplex}

Ist $I\su S$ von der Form
$I=S\setminus\{s_i\}$, so setzen wir
\begin{equation} Y_I:= \bigcup^\cdot_{h \in \B'^i} Y_{V^i}(w_h) .
\end{equation}
Da die Menge $\B'$ abgeschlossen ist, erh"alt man wiederum aus der
Proposition 3.3 und der Abschlu"seigenschaft der verallgemeinerten
Schubertzellen die nachstehende Folgerung.
\begin{Korollar}
Die Menge $Y_I$ ist abgeschlossen in $\Fl_g.$
\end{Korollar}
F"ur eine beliebige Teilmenge $I \subn S$ definieren wir
\begin{equation} Y_I:= \bigcap_{I \su J \atop \#(S\setminus J)=1} Y_J.
\end{equation}
Aus technischen Gr"unden setzen wir $Y_S=Y(\B').$ Diese Zuordnung liefert
also abgeschlossene Teilr"aume in $Y(\B')$, welche wir mit der reduzierten
Schemastruktur versehen. Ferner haben wir eine Operation der parabolischen
Untergruppe $P_I$ auf $Y_I,$ welche durch die nat"urliche  Aktion von $G$
auf $\Fl_g$ via Restriktion induziert wird.  Analog zu dem Fall der
$\Omega_I$  ist die obige
Zuordnung ebenfalls  inklusionserhaltend, d.h. aus $I \subset J$
folgt $Y_I \subset Y_J.$ F"ur zwei Teilmengen 
$I,J \subset S$ gilt sogar die Gleichheit 
\begin{equation}Y_{I \cap J}=Y_I \cap Y_J.
\end{equation}
F"ur ein beliebiges $I \subset S$ hat man zwischen $\Omega_I$ und $Y_I$
unter der Identifikation $\Fl_g = G/P_\mu$ die nachstehende Beziehung.
\begin{Lemma} 
Es gilt $Y_I= \bigcup\limits^\cdot_{w \in \Omega_I} BwP_\mu / P_\mu.$
\end{Lemma}
\beweis Sei zun"achst $I$ von der Gestalt $I=S \setminus \{s_i\}.$ Dann gilt
$$ Y_I=\bigcup^\cdot_{h \in \B'^i}Y_{V^i}(w_h) =\bigcup^\cdot_{h \in
  \B'^i}P_{\hat{i}}w_hP_\mu/P_\mu = \bigcup^\cdot_{h \in
  \B'^i}\bigcup_{w\in W_{\hat{i}}\dot{w}_h} BwP_\mu/P_\mu $$
$$ = \bigcup^\cdot_{h \in \B'^i}\bigcup^\cdot_{w \in W_{\hat{i}}\dot{w}_h \cap W^\mu} BwP_\mu/P_\mu =
      \bigcup^{\cdot}_{w \in\Omega_I} BwP_\mu/P_\mu.$$
Der allgemeine Fall ergibt sich dann aus der Vertr"aglichkeitseigenschaft
(6), (10) beider Zuordnungen bzgl. des Durchschnittes von Teilmengen in S.
\qed

Seien nun   $I \subset S$ und $g\in G(k)$ beliebig. Wir erhalten eine abgeschlossene Immersion
$$\phi_{g,I}: gY_I \hookrightarrow Y(\B') $$ von Variet"aten, wobei $gY_I$ das
Bild von $Y_I$ unter dem Translationsmorhpismus 
\Abb{g:}{Y(\B')}{Y(\B')}{x}{gx} sei. Da die parabolische
Untergruppe $P_I$ die Variet"at $Y_I$ fixiert, sind die Notationen $gY_I$
und $\Phi_{g,I}$ auch f"ur
$g\in (G/P_I)(k)$ sinnvoll.

Sei nun $F$ eine \'etale Garbe auf $Y(\B')$ und seien $I \subset J$ zwei
Teilmengen von $S$ mit $\#(J\setminus I)=1.$ Seien ferner $g \in (G/P_I)(k), h\in (G/P_J)(k)$ zwei
Linksnebenklassen, so da"s $g$ unter der kanonischen Projektion
$(G/P_I)(k) \longrightarrow (G/P_J)(k)$ auf $h$ abgebildet wird. Dann
sei $$p_{I,J}^{g,h}:  (\phi_{h,J})_*(\phi_{h,J})^* F \longrightarrow 
(\phi_{g,I})_*(\phi_{g,I})^* F$$ derjenige Morphismus von \'etalen Garben auf
$Y(\B'),$ der von der abgeschlossenen Immersion 
$$ gY_I\hookrightarrow hY_J$$
herr"uhrt. Wird $g$ nicht auf $h$ abgebildet, so setzen wir 
$p_{I,J}^{g,h}=0.$ Schlie"slich definieren wir $$ p_{I,J}= \!\!\!\!\!\!\!\!
\bigoplus_{(g,h) \in (G/P_I)(k)\times (G/P_J)(k)}\!\!\!\!\! p_{I,J}^{g,h}: \bigoplus_{h\in (G/P_J)(k)}\!\!\!\!\!(\phi_{h,J})_*(\phi_{h,J})^* F \longrightarrow \!\!\!\!\!
\bigoplus_{g \in (G/P_I)(k)}\!\!\!\!\! (\phi_{g,I})_*(\phi_{g,I})^* F.$$ F"ur
beliebige Teilmengen $I,J \subset S$ mit $\#J - \#I=1$ setzen wir
\[ d_{I,J} =\left\{ \begin{array}{l@{\quad : \quad }r}  (-1)^i p_{I,J} &
 J= I \cup \{s_i\} \\  0 & I \not\subset J \end{array} \right. .  \]
Wir erhalten einen Komplex 
$$(*): 0 \rightarrow F \rightarrow \bigoplus_{I \subset S \atop
\#(S\setminus I)=1} \bigoplus_{g \in (G/P_I)(k)} (\phi_{g,I})_*(\phi_{g,I})^* F \rightarrow
\bigoplus_{I \subset S \atop \#(S\setminus I)=2} \bigoplus_{g \in (G/P_I)(k)}
(\phi_{g,I})_*(\phi_{g,I})^* F \rightarrow $$
$$ \dots \rightarrow \bigoplus_{I \subset S \atop \#(S\setminus I)=d-2} \bigoplus_{g \in (G/P_I)(k)} (\phi_{g,I})_*(\phi_{g,I})^* F \rightarrow
\bigoplus_{g \in (G/B)(k)} (\phi_{g,\emptyset})_*(\phi_{g,\emptyset})^* F \rightarrow 0$$
von \'etalen Garben auf $Y(\B')$, in dem die Differentiale  durch die $d_{I,J}$
induziert werden.

Im n"achsten Paragraphen werden wir den folgenden Satz zeigen.
\begin{Satz} Sei $\B' \su \B$ eine abgeschlossene
Teilmenge mit $\B' \su \B^{ss}.$ Dann ist der Komplex $(*)$ azyklisch.
\end{Satz}

\section{Der Beweis von Satz 5.3 }

Die Azyklizit"at von $(*)$ wird mittels eines
Kontrahierbarkeitskriteriums f"ur simpliziale Komplexe, welche von partiell
geordneten Mengen definiert werden, gezeigt. 

Sei $X=(X,\leq)$ eine partiell geordnete Menge. Wir assoziieren zu $X$ einen
simplizialen Komplex 
$$ X^\bullet = \bigcup_{n \in \mathbb N} X^n,$$
wobei ein $n-$Simplex $\tau \in X^n$ durch ein geordnetes
$n-$Tupel
$$\tau=(x_0 < x_1 < \cdots < x_n)$$
mit Elementen $x_i \in X, \; i=0,\ldots, n$ gegeben ist. Insbesondere gilt
dann f"ur die 0-Simplizes $X^0$ von $X^\bullet$ die Identit"at $X^0=X.$

Ein Morphismus $f: X \longrightarrow Y$ von partiell
geordneten
Mengen ist eine Abbildung, welche ordnungserhaltend ist. Ein solcher
Morphismus
induziert eine simpliziale Abbildung
$$f^\bullet :X^\bullet \longrightarrow Y^\bullet$$
von simplizialen Komplexen. Durch diese Zuordnungen erh"alt man also einen
Funktor von der Kategorie der partiell geordneten Mengen in die Kategorie
der simplizialen Komplexe.
\begin{Beispiel}{\rm  Sei $V$ ein endlich-dimensionaler Vektorraum
"uber einem beliebigen K"orper $k$. Bezeichne mit $T_k(V)$ die Menge seiner
nicht-trivialen Unterr"aume. Man hat eine kanonische Ordnungsstruktur auf
$T_k(V)$, welche durch die Inklusion von Unterr"aumen gegeben ist. Ein
$n-$Simplex $\tau \in T_k(V)^n$ wird also durch eine Flagge
$$\tau =(\,(0) \subn U_0 \subn U_1 \subn \cdots \subn U_n \subn V \,)$$
von Unterr"aumen repr"asentiert. Dieser simpliziale Komplex ist bekanntlich
kanonisch isomorph zum  Tits-Komplex der Gruppe $GL(V)$ (vgl. \cite{Q} \S 3).}
\end{Beispiel}

F"ur simpliziale Komplexe, die durch partiell geordnete Mengen induziert
werden, hat man folgendes Kontrahierbarkeitskriterium von Quillen (vgl. \cite{Q}, 1.5.).
\begin{Proposition} Sei $X$ eine partiell geordnete Menge und $x_0
\in X$ ein \\ fixiertes Element. Falls ein Endomorphismus  $f: X \longrightarrow X$ existiert mit $$x\leq f(x) \geq x_0 \;\; \forall x \in X,$$ so ist der
zu $X$ assoziierte simpliziale Komplex $X^\bullet$ kontrahierbar.
\end{Proposition}

Wir wollen nun den Satz 5.3 beweisen.

\noindent {\bf Beweis von Satz 5.3} :  
Um die Behauptung zu beweisen, "uberpr"ufen wir die Exaktheit von
$(*)$ nach "Ubergang zur Lokalisierung in einem geometrischen Punkt 
$x \in Y(\B')(k^{sep})=Y(\B')(\overline{k}).$
Die Lokalisierung von $(*)$ in $x$ ergibt:
$$(**): 0 \longrightarrow F_x \longrightarrow \bigoplus_{I \subset S \atop
\#(S\setminus I)=1} \bigoplus_{g \in (G/P_I)(k) \atop
x \in gY_I(\overline{k})} F_x \longrightarrow
\bigoplus_{I \subset S \atop \#(S\setminus I)=2} \bigoplus_{g \in (G/P_I)(k)
  \atop x \in gY_I(\overline{k})} F_x \longrightarrow $$
$$ \dots \longrightarrow \bigoplus_{I \subset S \atop \#(S\setminus I)=d-2}
\bigoplus_{g \in (G/P_I)(k) \atop x \in gY_I(\overline{k})} F_x \longrightarrow
\bigoplus_{g \in (G/P_\emptyset)(k) \atop
x \in gY_\emptyset(\overline{k})} F_x \longrightarrow 0. $$
Wir definieren f"ur jedes $I \subn S$ die Menge
$$ R_I^x:=\{ gP_I(k)\in (G/P_I)(k); x \in gY_I(\overline{k})\} \su
(G/P_I)(k).$$ Sind  $I\su J \subn S$ zwei Teilmengen, so ist  $Y_I$ eine abgeschlossene
Untervariet"at von $Y_J$ (vgl. (10)).
Folglich gilt die Implikation
$$ gP_I(k)  \in R^x_I \Rightarrow gP_J(k) \in R^x_J,$$
d.h. mit jedem Simplex $\tau \in R^x_I$ vom Typ $I$ ist auch jedes
Untersimplex vom Typ $J$ in $R^x_J$ enthalten.
Dadurch erhalten wir einen simplizialen Unterkomplex des Tits-Komplexes
zur endlichen Gruppe $G(k)$, welchen wir mit $R^\bullet$ bezeichnen. 
Dabei repr"asentieren die Elemente aus
\begin{eqnarray*}
R^0 &=& \bigcup_{\{I; \#(S\setminus I)=1\}} R_I^x \;\;\;\mbox{ die 0-Simplizes},\\
R^1 &=& \bigcup_{\{I; \#(S\setminus I)=2\}} R_I^x \;\;\;\mbox{ die 1-Simplizes},\\
 & \vdots & \\
R^{d-2} &=& \bigcup_{\{I; \#(S\setminus I)=d-1\}} R_I^x  =
R^x_\emptyset\;\;\; \mbox{ die $(d-2)$-Simplizes}. 
\end{eqnarray*}

\noindent Im nachfolgenden Lemma wird gezeigt, da"s der simpliziale Komplex
$R^\bullet$ kontrahierbar ist. Daraus folgt die Azyklizit"at des zugeh"origen
Kettenkomplexes $(**)$ mit Werten in $F_x.$ Weil $x$ beliebig ist, folgt die Behauptung.
\qed

\begin{Lemma} Der simpliziale Komplex $R^\bullet$ ist kontrahierbar.
\end{Lemma}
\beweis Wir identifizieren 
den Tits-Komplex zu $G(k)$ mit dem simplizialen Komplex $T_k(V)^\bullet$
(vgl. Beispiel 6.1). Mittels dieser Identifikation lassen sich die Simplizes von $R^\bullet$
wie folgt beschreiben:
\begin{eqnarray*}
R^0 &=&\{ U \in T_k(V); x\in Y_U(h) \mbox{ wobei } h \in \B'^{\dim U}
\},\\ \\
R^1 &=&\{((0) \subn U_0 \subn U_1 \subn V) \in T_k(V)^1 ; x \in Y_{U_i}(h_i)
\mbox{ wobei }  h_i \in \B'^{\dim U_i}, \; i=0,1 \},\\
& \vdots & \\
R^{d-2} &=&\{((0) \subn U_0 \subn U_1 \subn \cdots \subn U_{d-2} \subn V)
\in T_k(V)^{d-2}; x \in Y_{U_i}(h_i)) \\
& & \mbox{ wobei }  h_i \in \B'^{\dim U_i},  \; i=0,\ldots,d-2 \}.
\end{eqnarray*}
Also ist $R^\bullet$ der zur partiell geordneten
Untermenge $R^0 \subset T_k(V)$ zugeh"orige simpliziale Komplex.
Wir wollen nun die Proposition  6.1 auf $R^0$
anwenden. Sei dazu 
$$\overline{T_k(V)}:= T_k(V) \cup \{(0), V\}.$$ Wir versehen
$\overline{T_k(V)}$ ebenfalls mit der offensichtlichen partiellen Ordnungsstruktur, so da"s $T_k(V)$ eine partiell geordnete Untermenge von
$\overline{T_k(V)}$ ist.
Sei nun $U_0 \in R^0$ ein minimales Element, d.h. $U_0$ enth"alt au"ser
sich selbst kein anderes Objekt aus $R^0.$ Betrachte die
Abbildung     
\Abb{f:}{R^0}{\overline{T_k(V)}\;\;.}{U}{U_0+U}
Diese Abbildung stellt offenbar einen Morphismus von partiell geordneten
Mengen dar, und es gilt
$$ U \leq f(U) \geq U_0, \;\; \forall \,U \in R^0     .$$
\noindent Nach der Proposition 6.2 reicht es also zu zeigen, da"s das Bild
von $f$ in $R^0$ enthalten ist. Sei also $U \in R^0,$ d.h f"ur die
zugeh"orige Unterfunktion $h \in \B^{\dim U}$ mit $x\in Y_U(h)$ gilt $h \in
\B'$ (vgl. (2)). Sei $h'$ die entsprechende Unterfunktion der L"ange $\dim(f(U))$
mit $x \in Y_{f(U)}(h').$ Wegen der Abgeschlossenheit von $\B'$
gen"ugt es die Ungleichung
\begin{equation} \deg_x(f(U)) = \sum_x xh'(x) \geq \sum_x xh(x) =
  \deg_x(U)
\end{equation}
zu zeigen. Beachte, da"s wegen der Inklusion $\B' \subset \B^{ss}$ und der
Identit"at $\deg_x(V) =0$ der Fall $f(U)=V$ dann nicht m"oglich ist. 
Um die Ungleichung (11) zu zeigen, hat man die folgenden drei F"alle zu betrachten.

\medskip
{\noindent \bf 1. Fall:} Es ist $U_0 \su U.$ Dann gilt
$$ \deg_x(f(U))=\deg_x(U) .$$

{\noindent \bf 2. Fall:} Es ist $U_0 \cap U = (0).$ Dann gilt nach Korollar 1.10
$$ \deg_x(f(U))=\deg_x(U+U_0) \geq \deg_x(U) + \deg_x(U_0)  > 0 .$$

{\noindent \bf 3. Fall:} Es ist $U_0 \cap  U \neq (0)$ und $U_0 \nsu U
$. Wegen der speziellen Wahl von $U_0$ hat man dann
$$ \deg_x(U \cap U_0) \leq 0 .$$
Somit gilt wiederum nach Korollar 1.10
$$ \deg_x(f(U))=\deg_x(U+U_0) \geq \deg_x(U) + \deg_x(U_0) - \deg_x(U \cap U_0)
>0.\mbox{          \qed}$$ 

\section{Der Beweis der S"atze 4.1 und 4.2}

Sei $I\subn S$ eine Teilmenge.
F"ur eine Zahl $0 \leq i \leq m_I:=\max\{l(w); w \in \Omega_I\}$
definieren wir die Untermenge
$$\Omega^i_I=\{w \in \Omega_I; l(w)\leq i\}$$ von $\Omega_I.$
\begin{Proposition}
Es gilt
\begin{equation} H^\ast_{\acute{e}t}(Y_I)=\bigoplus_{w\in \Omega_I}
\mathbb Q_\ell(-l(w))[-2l(w)].
\end{equation}
\end{Proposition}
F"ur den Beweis ben"otigen wir folgendes Lemma.

\begin{Lemma} Sei $X$ eine Variet"at "uber $k$ und
$$ X=X^m\supset X^{m-1}\supset \cdots \supset X^1 \supset X^0\supset X^{-1}=\emptyset$$
eine absteigende Sequenz von abgeschlossenen Untervariet"aten "uber $k$. Sei f"ur jedes 
$0 \leq i \leq m$ 
$$ X^i\setminus X^{i-1}= \coprod_{J_i} \mathbb A_k^i$$ 
eine endliche direkte Summe von standard-affinen R"aumen der Dimension $i$
bzgl. einer Indexmenge $J_i$.
Dann gilt
$$H_c^\ast(X)=\bigoplus_{i=0}^m H^\ast_c(X^i\setminus X^{i-1})=\bigoplus_{i=0}^m
\bigoplus_{J_i} H^{2i}_c(\mathbb A_k^i)=\bigoplus_{i=0}^m \bigoplus_{J_i}\mathbb Q_\ell(-i)[-2i].$$
\end{Lemma}

\beweis Es gilt
$$ H^j_c(\mathbb A^i_k)=\left\{ \matrix{ \mathbb Q_\ell(-i) & j=2i \cr 0 & \mbox{sonst} }
\right. $$
  Da die Kohomologie von $\mathbb A_k^i$ also im Grad $2i$ konzentriert ist,
folgt die Behauptung unmittelbar aus der sukzessiven Anwendung der
langen exakten Kohomologiesequenzen zu den Tripeln
$$ X^i\setminus X^{i-1} \hookrightarrow X^i \hookleftarrow
X^{i-1}.\;\square$$

{\noindent \bf Beweis von Proposition 7.1:}  F"ur eine Zahl
$0 \leq i \leq m_I$ definieren wir
die folgenden Untermengen von $Y_I:$
$$ Y_I^i:= \bigcup_{w \in \Omega^i_I} BwP_\mu/P_\mu .$$
Offenbar ist diese Untervariet"at der Durchschnitt von $Y_I$ und der Variet"at
$\bigcup_{w\in W \atop l(w)\leq i} BwP_\mu/P_\mu.$ Letztere ist aber
abgeschlossen in $\Fl_g.$ Also handelt es sich bei den
$Y_I^i$ um abgeschlossene Untervariet"aten der $Y_I.$
Wir erhalten somit eine absteigende Sequenz von abgeschlossenen Untervariet"aten
$$ Y_I=Y_I^{m_I} \supset Y_I^{m_I-1} \supset \cdots \supset Y_I^0 \supset Y_I^{-1}:=\emptyset.$$ 
Wegen der Abschlu"seigenschaft der Bruhatzellen bzgl. der Bruhatordnung
zerlegen sich  die $Y_I^i\setminus
Y_I^{i-1} $ in die direkte Summe ihrer  Bruhatzellen, d.h.
$$ Y_I^i\setminus Y_I^{i-1} =\coprod_{w \in \Omega_I \atop l(w)=i}
BwP_\mu/P_\mu \;, \; i=0,\ldots,m_I.$$
Nun sind die Bruhatzellen $BwP_\mu/P_\mu $ aber standard-affine R"aume der
Dimension $l(w)$. 
Die Anwendung von Lemma 7.2 liefert die Behauptung.
\qed

\medskip
Im folgenden soll f"ur ein Weylgruppenelement $w \in W^\mu$ und eine
Teilmenge $I \subset S$ mit $H(Y_I,w)$ der Beitrag von $w$ zur direkten
Summe (12) bezeichnet werden. Es ist also
\begin{equation}
H(Y_I,w)=\left\{ 
\begin{array}{r@{\quad:\quad}l} \mathbb Q_\ell(-l(w))[-2l(w)] & w \in \Omega_I \\ 
0 & \mbox{ sonst } \end{array}
\right. .
\end{equation}
Dann gilt offenbar
\begin{equation}
 H^\ast_{\acute{e}t}(Y_I)= \bigoplus_{w\in W^\mu} H(Y_I,w).
\end{equation}
Seien nun $I\subset J $ zwei echte Teilmengen von S. Dann
ist $Y_I$ abgeschlossen in $Y_J$. Wir betrachten den Restriktionshomomorphismus
$$\phi_{I,J}: H^{\ast}_{\acute{e}t}(Y_J) \longrightarrow
H^{\ast}_{\acute{e}t}(Y_I).$$
Wegen (14) erhalten wir eine Graduierung von $\phi_{I,J},$
$$\phi_{I,J}=\bigoplus_{(w,w') \in W^\mu \times W^\mu} \phi_{w,w'}:
\bigoplus_{w\in W^\mu} H(Y_J,w) \longrightarrow \bigoplus_{w'\in W^\mu}
H(Y_I,w').$$
Aus der Konstruktion in Lemma 7.2 ergeben sich f"ur die $\phi_{w,w'}$
folgende Werte:

\begin{equation} \phi_{w,w'}=\left\{\begin{array}{r@{\quad:\quad}l} id  & w=w'
    \\ 0 & w \neq w' \end{array} \right. . \end{equation}

Vor dem Beweis von Satz 4.1 ben"otigen wir noch ein Resultat von
Lehrer bzw. Bj"orner. Wir konstruieren einen Komplex, der analog zur Sequenz $(*)$ aufgebaut ist.
Seien  $I\subset J
\subset S$ zwei Teilmengen mit $\#(J\setminus I) =1.$ Wir erhalten 
dann einen Homomorphismus 
$$p_{I,J}:i^G_{P_J} \longrightarrow i^G_{P_I},$$ welcher durch
die Surjektion $(G/P_I)(k) \longrightarrow (G/P_J)(k)$ gegeben wird.
F"ur zwei beliebige Mengen $I,J \subset S$ mit $\#J -\#I=1$  definieren wir 
$$d_{I,J}=\left\{ \matrix{ (-1)^i p_{I,J} & J = I  \cup \{s_i\} \cr 0 & I \not\subset J } \right. .$$ 
Wir erhalten somit f"ur jedes $I_0 \subset S $ einen $\mathbb
Z$-indizierten Komplex
$$K_{I_0}^\bullet: 0 \rightarrow i^G_G \rightarrow
  \bigoplus_{I_0 \subset I \subset S \atop \#(S\setminus I)=1} i^G_{P_I}  \rightarrow
\bigoplus_{I_0 \subset I \subset S \atop \#(S\setminus I)=2} i^G_{P_I} \rightarrow
\dots \rightarrow \bigoplus_{I_0 \subset I \subset S \atop \#(S\setminus
  I)=\#(S\setminus I_0)-1 } i^G_{P_I}  \rightarrow i^G_{P_{I_0}},$$
in dem die Differentiale durch die  $d_{I,J}$ gegeben werden. Dabei steht
die Komponente $i^G_G$ im Grad $-1.$ Der Komplex $K_{I_0}^\bullet$
identifiziert sich, bis auf den Term $i^G_G$, gerade mit demjenigen
Kettenkomplex mit Werten in $\mathbb Q_\ell,$
welcher durch den simplizialen Komplex $\Delta_{I_0}$ aus \cite{L} 
bzw. \cite{Bj} induziert wird.

\begin{Satz}{\rm (\cite{L} 4.3 bzw. \cite{Bj} 4.2 )} Der Komplex $K_{I_0}^\bullet$ ist azyklisch.
\end{Satz}   
Wir erw"ahnen auch folgendes wohlbekanntes Lemma.
\begin{Lemma} Jede Erweiterung des $Gal(\overline{k}/k)$-Moduls $\mathbb
  Q_\ell(m)$ durch $\mathbb Q_\ell(n)$ mit $m\neq n$ zerf"allt.
\end{Lemma}


Nun haben wir alle Hilfsmittel beisammen, um Satz 4.1 zu beweisen.

{\noindent \bf Beweis von Satz 4.1:} Wir setzen
$$\overline{Y(\B')}:=Y(\B')\times _k \overline{k} \mbox{ bzw. }
\overline{\Phi_{g,I}}=\Phi_{g,I} \times _k id_{\overline{k}}.$$ Wir betrachten
die aus der exakten Sequenz $(*)$ resultierende Spektralsequenz
$$ E_1^{p,q} = H^q_{\acute{e}t}(\overline{Y(\B')},\bigoplus_{I \subset S \atop
\#(S\setminus I)=p+1 } \bigoplus_{g \in (G/P_I)(k)}(\overline{\phi_{g,I}})_*(\overline{\phi_{g,I}})^* \mathbb Q_\ell) \Longrightarrow
H^{p+q}_{\acute{e}t}(\overline{Y(\B')},\mathbb Q_\ell).$$ 
Es gilt
\begin{eqnarray*} E_1^{p,q} &=&
  H^q_{\acute{e}t}(\overline{Y(\B')},\bigoplus_{I \subset S \atop \#(S\setminus I)=p+1 } \bigoplus_{g \in
(G/P_I)(k)}(\overline{\phi_{g,I}})_*(\overline{\phi_{g,I}})^* \mathbb Q_\ell) \\
&=& \bigoplus_{I \subset S \atop \#(S\setminus I)=p+1 } \bigoplus_{g \in
  (G/P_I)(k)}H^q_{\acute{e}t}(\overline{Y_I},(\overline{\phi_{g,I}})^* \mathbb Q_\ell)
=\bigoplus_{I \subset S \atop \#(S\setminus I)=p+1 } \bigoplus_{(G/P_I)(k)}
H_{\acute{e}t}^q(Y_I).
\end{eqnarray*}
Die Anwendung von $(14)$ und $(15)$ liefert eine Zerlegung des
$E_1$-Terms der Spektralsequenz in Unterkomplexe 
$$ E_1=\bigoplus_{w \in W^\mu} E_{1,w}$$ mit 
$$E_{1,w}^{p,q}=\left\{ \mwmatrix{\bigoplus_{I \subset S \atop
      \#(S\setminus I)=p+1} \bigoplus_{(G/P_I)(k)} H(Y_I,w) & \;\;\;\;\; q=2l(w) \cr
    0 & \;\;\;\;\; q\neq 2l(w) } \right. .$$
Also ist $E_{1,w}$ der Unterkomplex
$$ E_{1,w}: \bigoplus_{I \subset S \atop \#(S\setminus I)=1} \bigoplus_{
(G/P_I)(k)}H(Y_I,w)\rightarrow \bigoplus_{I \subset S
\atop \#(S \setminus I)=2} \bigoplus_{(G/P_I)(k)}H(Y_I,w) \rightarrow
\dots \rightarrow\bigoplus_{(G/B)(k)}H(Y_\emptyset,w).$$
F"ur eine Teilmenge $I \subn S$ gilt nach (7) und (13)
\[ H(Y_I,w)=\left\{ 
\begin{array}{r@{\quad:\quad}l} \mathbb Q_\ell(-l(w))[-2l(w)] &  I_w \su I
\\ 0 & \mbox{sonst} \end{array}
\right. . \]
Deshalb vereinfacht sich $E_{1,w}$ zu 
$$ E_{1,w}: \Big(\bigoplus_{I_w \subset I \atop \#(S\setminus I)=1}i^G_{P_I}
\rightarrow \bigoplus_{I_w \subset I \atop \#(S\setminus I) =2}
i^G_{P_I} \rightarrow \dots \rightarrow i^G_{P_{I_w}}\Big)(-l(w))[-2l(w)].$$
Also stimmt $E_{1,w}$ bis auf den Term $i^G_G$, dem Tate-Twist $(-l(w))$
und dem Shift $[-2l(w)]$ "uberein mit dem Komplex $K_{I_w}^\bullet.$
Genauer gesagt haben wir die folgende exakte Sequenz von Komplexen:
$$ 0 \rightarrow i^G_G(-l(w))[-2l(w)+1] \rightarrow K_{I_w}^\bullet(-l(w))
[-2l(w)]\rightarrow E_{1,w} \rightarrow 0.$$

\noindent Somit ergeben sich f"ur den $E_{2,w}$-Term nachstehende drei F"alle:

\medskip
$\begin{array}{rclrl} I_w=S & : & E_{2,w}^{p,q} & =& 0 \;\;  p\geq 0, q\geq 0 \\ \\
\#(S \setminus I_w)=1 & : & E_{2,w}^{0,2l(w)} &=&  i^G_{P_{I_w}}(-l(w)) \\ \\
           &   & E_{2,w}^{p,q} &=&  0 \;\; (p,q)\neq(0,2l(w)) \\ \\
\#(S \setminus I_w) > 1 & : & E_{2,w}^{0,2l(w)} & = &  i^G_G(-l(w)) \\ \\
& & E_{2,w}^{j,2l(w)} &=&  0 \;\; j=1,\ldots,\#(S \setminus I_w)-2 \\ \\
& & E_{2,w}^{j,2l(w))} &=& v^G_{P_{I_w}}(-l(w)) \;\; j=\#(S\setminus
I_w)-1\\ \\
& & E_{2,w}^{p,q} &=& 0 \;\;  q\neq 2l(w) \mbox{ oder }  p>\#(S \setminus I_w)-1.
\end{array}$

\medskip
\noindent Nun haben die  $E_2^{p,q}\neq (0)$ den Tate-Twist $-\frac{q}{2}$. Da aber jeder Homomorphismus von
Galois-Moduln unterschiedlichen Tate-Twistes trivial ist, stellt der
$E_2$-Term bereits  den
$E_\infty$-Term dar. Zusammengefa"st gilt also f"ur jedes $n\geq 0$
$$gr^p(H_{\acute{e}t}^n(Y(\B'),\mathbb Q_\ell)=E_\infty^{p,n-p}=E_2^{p,n-p}=\bigoplus_{w\in W^\mu} E_{2,w}^{p,n-p}$$
\[ =\left\{ 
\begin{array}{r@{\quad:\quad}l}
  \bigoplus\limits_{w \in W^\mu, \#(S\setminus I_w)=1 \atop
    2l(w)=n}i^G_{P_{I_w}}(-l(w))\oplus \bigoplus\limits_{w \in W^\mu, \#(S
    \setminus I_w)>1 \atop 2l(w)=n} i^G_G(-l(w)) &  p=0 \\
\bigoplus\limits_{w \in W^\mu \atop 2l(w) + \#(S\setminus I_w)-1=n}
v^G_{P_{I_w}}(-l(w))  &
  p>0 \\
\end{array}
\right.  \] 
Nun sind aber nach Lemma 7.4 Erweiterungen von $\mathbb Q_\ell(m)$ durch
$\mathbb Q_\ell(n)$ mit $m\neq n$ trivial.
Somit gilt
$$ H_{\acute{e}t}^n(Y(\B'))\cong \bigoplus_{p \in
  \mathbb N} gr^p(H_{\acute{e}t}^n(Y(\B')))$$
\[=\bigoplus_{w\in W^\mu, \#(S\setminus I_w)=1 \atop
  2l(w)=n}\!\!\!\!\!i^G_{P_{I_w}}(-l(w))\oplus\!\!\!\!\!\! \bigoplus_{w\in W^\mu, \#(S\setminus
  I_w) >1 \atop 2l(w)=n}\!\!\!\!\!i^G_G(-l(w)) \oplus \!\!\!\!\!
\bigoplus_{w \in W^\mu \atop 2l(w) + \#(S\setminus I_w)-1=n}\!\!\!\!\! v^G_{P_{I_w}}(-l(w)) .\]
Hieraus folgt die Behauptung.
\qed

\bigskip
Indem wir das Ergebnis von Satz 4.1 verwenden, l"a"st sich nun Satz 4.2 beweisen.

{\noindent \bf Beweis von Satz 4.2:}
Es soll zun"achst der Homomorphismus
$$i^\ast: H_{\acute{e}t}^\ast(\Fl_g) \longrightarrow
H_{\acute{e}t}^\ast(Y(\B'))$$
bestimmt werden, der durch die Inklusion $Y(\B') \stackrel{i}{\hookrightarrow} \Fl_g$ induziert
wird. Dann folgt n"amlich, wie wir sehen werden, die Behauptung unmittelbar aus der zum Tripel
$$ \Fl_g(\B') \hookrightarrow \Fl_g \hookleftarrow Y(\B')$$
geh"origen langen exakten Kohomologiesequenz.
Betrachte den Doppelkomplex $(\tilde{E_1^\bullet},\partial^\bullet,\partial'^\bullet)$ definiert durch

$$\tilde{E_1}^{p,q}=\left\{ \matrix{ H_{\acute{e}t}^q(\Fl_g)= \bigoplus\limits_{w
      \in W^\mu \atop  2l(w)=q} \mathbb Q_\ell(-l(w))[-2l(w)] & p=0 \cr \cr 0 & p>0 }\right. ,$$
\[(\partial^{p,q}:\tilde{E_1}^{p,q} \longrightarrow \tilde{E_1}^{p+1,q}) =0
\;\; \forall (p,q) \in \mathbb N \times \mathbb N, \]  
$$(\partial'^{p,q}: \tilde{E_1}^{p,q} \longrightarrow  \tilde{E_1}^{p,q+1}) =0
\;\; \forall (p,q) \in \mathbb N \times \mathbb N . $$
Offenbar definiert $\tilde{E_1^\bullet}$ eine Spektralsequenz, welche gegen
$H_{\acute{e}t}^\ast(\Fl_g)$ konvergiert.
Sei f"ur jedes $I\subset S$ mit $\#(S\setminus I)=1$ 
$$D_I^\ast : H_{\acute{e}t}^\ast(\Fl_g) \longrightarrow
\bigoplus_{(G/P_I)(k)} H_{\acute{e}t}^\ast(Y_I)$$ der Diagonalhomomorphismus, der durch die
abgeschlossene Immersion $Y_I \hookrightarrow \Fl_g$ induziert wird. Betrachte den Morphismus von Doppelkomplexen
$$f_1^\bullet : \tilde{E_1^\bullet} \longrightarrow E_1^\bullet $$
definiert durch
\[ f_1^{p,q} =\left\{ 
\begin{array}{r@{\quad:\quad}l}
\oplus_{I} D_I^q & p=0  \\
 0 & p>0. \end{array}
\right. \] 
Bezeichne mit $s(E_1)^{\bullet}$ bzw. $s(\tilde{E_1})^\bullet$ die
assoziierten einfachen Totalkomplexe und mit
$$s(f^\bullet_1): s(\tilde{E_1})^\bullet \longrightarrow s(E_1)^\bullet$$  den durch $f^\bullet_1$
induzierten Homomorphismus. Bezeichne mit  $\mbox{cone}(s(f^\bullet_1))^\bullet$ den Kegel von $s(f_1^\bullet).$
Offenbar identifiziert sich $s(f_1^\bullet)$ in der Kohomologie gerade 
mit dem Homomorphismus 
$$i^\ast: H_{\acute{e}t}^\ast(\Fl_g) \longrightarrow
H_{\acute{e}t}^\ast(Y(\B')).$$
Deshalb berechnet die Kohomologie des geshifteten Kegels $\mbox{cone}(s(f^\bullet_1))^\bullet[-1]$
die gew"unschten Kohomologiegruppen $H^n_c(\Fl_g(\B')), n \geq 0.$
Wie man sofort verifiziert gilt aber
$$ \mbox{cone}(s(f^\bullet_1))^\bullet = \bigoplus_{w \in W^\mu}
K_{I_w}^\bullet(-l(w))[-2l(w)] .$$ 
Unter Ausnutzung der Gleichheit $K_{I_w}^{-1}=i^G_G$  folgt die Behauptung.\qed

\vspace{1cm}

\centerline{Mathematisches Institut der Universit"at zu K"oln} 
\centerline{Weyertal   86-90, 50931 K"oln, Deutschland}
\centerline{e-mail: sorlik@mi.uni-koeln.de}

\end{document}